%&amstex          
\input amstex\documentstyle{amsppt}  
\pagewidth{12.5cm}\pageheight{19cm}\magnification\magstep1
\topmatter
\title On certain varieties attached to a Weyl group element\endtitle
\author G. Lusztig\endauthor
\address{Department of Mathematics, M.I.T., Cambridge, MA 02139}\endaddress
\thanks{Supported in part by the National Science Foundation}\endthanks
\endtopmatter   
\document

\define\we{\wedge}

 \define\db{\dot b} \define\dc{\dot c} \define\dv{\dot v}

\define\dw{\dot w}

\define\ds{\dot s}

\define\bcp{\bar{\cp}}

\define\hs{\hat s}
\define\hg{\hat g}
\define\hw{\hat w}

\define\pe{\perp}
\define\si{\sim}

\define\sqc{\sqcup}
\define\ovs{\overset}
\define\qua{\quad}

\define\hG{\hat G}

\define\tis{\ti\s}

\define\lb{\linebreak}

\define\op{\oplus}

\define\part{\partial}

\define\n{\notin}

\define\m{\mapsto}
\define\do{\dots}

\define\bsl{\backslash}

\define\lra{\leftrightarrow}

\define\sub{\subset}    

\define\T{\times}
\define\ti{\tilde}
\define\nl{\newline}
\redefine\i{^{-1}}

\define\ov{\overline}

\define\bbq{\bar{\QQ}_l}

\define\Ad{\text{\rm Ad}}

\define\a{\alpha}
\redefine\b{\beta}

\define\g{\gamma}
\redefine\d{\delta}
\define\e{\epsilon}

\redefine\o{\omega}
\define\p{\pi}

\define\r{\rho}
\define\s{\sigma}
\redefine\t{\tau}

\define\k{\kappa}
\redefine\l{\lambda}
\define\z{\zeta}
\define\x{\xi}

\redefine\G{\Gamma}

\redefine\L{\Lambda}
\define\Ph{\Phi}
\define\Ps{\Psi}

\define\ii{\bold i}

\define\kk{\bold k}

\define\nn{\bold n}

\define\qq{\bold q}

\redefine\ss{\bold s}

\define\NN{\bold N}

\define\QQ{\bold Q}

\define\UU{\bold U}

\define\WW{\bold W}
\define\ZZ{\bold Z}
\define\XX{\bold X}

\define\cb{\Cal B}

\define\cf{\Cal F}
\define\cg{\Cal G}
\define\ch{\Cal H}
\define\ci{\Cal I}

\define\cl{\Cal L}

\define\co{\Cal O}
\define\cp{\Cal P}

\define\car{\Cal R}

\define\ct{\Cal T}

\define\fB{\frak B}

\define\fF{\frak F}

\define\fZ{\frak Z}

\define\tz{\ti z}

\define\tT{\ti T}

\define\tX{\ti X}

\define\tZ{\ti Z}

\define\bul{\bullet}

\define\tfB{\ti\fB}    
\define\tXX{\ti{\XX}}

\define\BR{BR}
\define\MB{MB}
\define\BM{BM}
\define\DL{DL}
\define\DM{DM}
\define\GM{GM}
\define\GP{GP}
\define\GKP{GKP}
\define\HE{H1}
\define\HEA{H2}
\define\COX{L1}
\define\RFC{L2}
\define\CS{L3}
\define\CSD{L4}
\define\WEU{L5}
\define\OR{OR}

\head 0. Introduction and statement of results\endhead
\subhead 0.1\endsubhead
Let $\kk$ be an algebraically closed field. Let $G$ be a connected reductive algebraic group over $\kk$. We 
assume that we are in one of the following two cases.

(1): $G$ is the identity component of a reductive group $\hG$ with a fixed connected component $D$.

(2): $\kk$ is an algebraic closure of a finite field 
$F_q$ and $G$ has a fixed $F_q$-rational structure with Frobenius map $F:G@>>>G$.
\nl
In case (1) we set $q=1$ and denote by $F:G@>>>G$ the identity map of $G$ so that $G^F=G$. Thus when $q=1$ we are
in case (1) and when $q>1$ we are in case (2).

Let $\cb$ be the variety of Borel subgroups of $G$. Let $\WW$ be an indexing set for the set of $G$-orbits on 
$\cb\T\cb$ for the diagonal $G$-action. Let $\co_w$ be the $G$-orbit corresponding to $w\in\WW$. Note that $\WW$ 
is naturally a Coxeter group with length function $l(w)=\dim\co_w-\dim\cb$. 

Let $I$ be an indexing set for the set $S$ of simple reflections of $\WW$. Let $s_i\in S$ be the simple 
reflection corresponding to $i\in I$. For $B\in\cb$ we have $gBg\i\in\cb$ for any $g\in D$ (if $q=1$) and 
$F(B)\in\cb$ (if $q>1$). There is a unique automorphism of $\WW$ (denoted by $\bul$ or by $w\m w^\bul$) such that
$\co_{w^\bul}=g\co_wg\i$ for all $w\in\WW,g\in D$ (if $q=1$) and $\co_{w^\bul}=F(\co_w)$ for all $w\in\WW$ (if 
$q>1$). We have $l(w^\bul)=l(w)$ for all $w\in\WW$. Hence there is a unique bijection $i\m i^\bul$ of $I$ such 
that $s_i^\bul=s_{i^\bul}$ for all $i\in I$.

Two elements $w,w'\in\WW$ are said to be $\bul$-conjugate if $w'=a\i wa^\bul$ for some $a\in\WW$. The relation of
$\bul$-conjugacy is an equivalence relation on $\WW$; the equivalence classes are said to be $\bul$-conjugacy 
classes. A $\bul$-conjugacy class $C$ in $\WW$ (or an element of it) is said to be $\bul$-elliptic if $C$ does 
not meet any $\bul$-stable proper parabolic subgroup of $\WW$ (see \cite{\HE}). (In the case where $\bul=1$ we 
say "elliptic, conjugacy class" instead of "$\bul$-elliptic, $\bul$-conjugacy class".) For $w\in\WW$ let 

$\fB_w=\{(g,B)\in D\T\cb;(B,gBg\i)\in\co_w\}$ (if $q=1$)

$X_w=\{B\in\cb;(B,F(B))\in\co_w\}$ (if $q>1$).
\nl
This is naturally an algebraic variety over $\kk$. (The variety $X_w$ is defined in \cite{\DL}. The variety 
$\fB_w$ appears in \cite{\CS} assuming that $D=G$ and in \cite{\CSD} in general.) We shall use the notation
$\XX_w$ for either $\fB_w$ or $X_w$. Let $\r:\fB_w@>>>D$ be the first projection.

Now $G^F$ acts on $\XX_w$ by $x:(g,B)\m(xgx\i,xBx\i)$ (if $q=1$) and by $x:B\m xBx\i$ (if $q>1$).

One of the themes of this paper is the analogy between $X_w$ and $\fB_w$. It seems that $\fB_w$ is a limit case 
of $X_w$ as $q\to1$. For example it is likely that for any $i$, the multiplicities of various unipotent character
sheaves on $D$ in a Jordan-H\"older series of the $(i+\dim G)$-th perverse cohomology sheaf of $\r_!\bbq$ (with 
$q=1$) are the same as the multiplicities of various irreducible unipotent representations of $G^F$ in the 
$G^F$-module $H^i_c(X_w,\bbq)$ (with $q>1$). Here $l$ is a fixed prime number invertible in $\kk$.

\subhead 0.2\endsubhead
From \cite{\DL, 1.11} it is known that if $w\in\WW$, $X_w$ has a natural finite covering $\tX_w$. We now show 
that (at least if $w$ is $\bul$-elliptic and $G$ is semisimple), $\fB_w$ has a natural finite covering $\tfB_w$.

Let $B^*\in\cb$ and let $T^*$ be a maximal torus of $B^*$; if $q>1$ we assume in addition that $B^*,T^*$ are 
defined over $F_q$. Let $U^*$ be the unipotent radical of $B^*$. If $q=1$ let $d\in D$ be such that 
$dT^*d\i=T^*,dB^*d\i=B^*$. Let $N=\{n\in G;nT^*n\i=T^*\}$. We identify $N/T^*=\WW$ by $nT^*\lra w$, 
$(B^*,nB^*n\i)\in\co_w$. According to Tits, for each $w\in\WW$ we can choose a representative $\dw\in N$ in such 
a way that $\dw=\dw_1\dw_2$ whenever $w,w_1,w_2$ in $\WW$ satisfy $w=w_1w_2$, $l(w)=l(w_1)+l(w_2)$. We can also 
assume that, if $w'=w^\bul$, then $\dw'=d\dw d\i$ (if $q=1$) and $\dw'=F(\dw)$ (if $q>1$).
For $w\in\WW$ let $U^*_w=U^*\cap\dw U^*\dw\i$ and let $T^*_w=\{t_1\in T^*;\dw\i t\dw=dtd\i\}$ (if $q=1$), 
$T^*_w=\{t\in T^*;\dw\i t\dw=F(t)\}$ (if $q>1$); let 
$$\tfB_w=\{(g,g'U^*_w)\in D\T G/U^*_w;g'{}\i gg'\in\dw U^*d\}\text{ (if }q=1),$$
$$\tX_w=\{g'U^*_w\in G/U^*_w;g'{}\i F(g')\in\dw U^*\}\text{ (if }q>1).$$
We shall use the notation $\tXX_w$ for either $\tfB_w$ or $\tX_w$. Now $G^F$ acts on $\tXX_w$ by 
$x:(g,g'U^*_w)\m(xgx\i,xg'U^*_w)$ (if $q=1$) and by $x:g'U^*_w\m xg'U^*_w$ (if $q>1$). Also $T^*_w$ acts (freely)
on $\tXX_w$ by $t:(g,g'U^*_w)\m(g,g't\i U^*_w)$ (if $q=1$) and by $t:g'U^*_w\m g't\i U^*_w$ (if $q>1$); this 
action commutes with the $G^F$-action. Define $\p_w:\tXX_w@>>>\XX_w$ by $(g,g'U^*_w)\m(g,g'B^*g'{}\i)$ (if $q=1$)
and by $g'U^*_w\m g'B^*g'{}\i$ (if $q>1$). Note that $\p_w$ is compatible with the $T^*_w$ action where $T^*_w$ 
acts on $\XX_w$ trivially. Let $\frak F$ be the fibre of $\p_w$ at a point of $\XX_w$. Then for some 
$t_0\in T^*$, $\frak F$ can be identified with $\{t\in T^*;\Ad(d)(t\i)\Ad(\dw\i)(t)=t_0\}$ (if $q=1$) and with 
$\{t\in T^*;F(t)\i\Ad(\dw\i)(t)=t_0\}$ (if $q>1$); hence it is either empty or a principal homogeneous space for
$T^*_w$. Now if $q>1$, $T^*_w$ is finite, hence the homomorphism $T^*@>>>T^*$, $t\m F(t)\i\Ad(\dw\i)(t)$ is 
surjective and $\frak F$ is a principal homogeneous space for $T^*_w$ so that in this case, $\p_w$ is a principal
$T^*_w$-bundle. If for $q=1$ we assume that $G$ is semisimple and $w$ is $\bul$-elliptic then $T^*_w$ is finite,
hence the homomorphism $T^*@>>>T^*$, $t\m\Ad(d)(t\i)\Ad(\dw\i)(t)$ is surjective and $\frak F$ is a principal 
homogeneous space for $T^*_w$ so that in this case, $\p_w$ is again a principal $T^*_w$-bundle. 

Here is another reason why $\fB_w$ looks like $X_w$ when $q\to1$: assuming that $w$ is $\bul$-elliptic and $G$ is
semisimple, the number $|T^*_w|$ (in the case $q=1$) is obtained from the number $|T^*_w|$ (in the case $q>1$) 
viewed as a polynomial in $q$ by substituting $q=1$.

The following result gives another instance of analogous behaviour of $X_w,\fB_w$.

\proclaim{Theorem 0.3} Assume that $w\in\WW$ is $\bul$-elliptic and that $w$ has minimal length in its 
$\bul$-conjugacy class. If $q=1$ assume further that $G$ is semisimple. 

(a) If $q=1$ (resp. $q>1$), any isotropy group of the $G^F$ action on $\tfB_w$ (resp. $\tX_w$) is $\{1\}$.

(b) If $q=1$ (resp. $q>1$), any isotropy group of the $G^F$ action on $\fB_w$ (resp. $X_w$) is isomorphic to a 
subgroup of $T^*_w$; hence it is a finite diagonalizable group.

(c) If $q=1$, the varieties $\fB_w$ and $\tfB_w$ are affine.
\endproclaim
Note that (c) has the following known analogue (see \cite{\DL} for sufficiently large $q$ and 
\cite{\OR}, \cite{\HE} for any $q$): 

(d) {\it If $q>1$, the varieties  $X_w$ and $\tX_w$ are affine.}
\nl
The proof of the theorem (given in \S3) extends the proof of a weaker form of (b) given in \cite{\WEU, 5.2}.

Let $G\bsl\tfB_w$ (resp. $G^F\bsl\tX_w$) be the set of orbits of the $G$-action on $\tfB_w$ (resp. of the 
$G^F$-action on $\tX_w$). Let $G\bsl\fB_w$ (resp. $G^F\bsl X_w$) be the set of orbits of the $G$-action on 
$\fB_w$ (resp. of the $G^F$-action on $X_w$). By (a)-(c) above. $G\bsl\tfB_w$ and $G\bsl\fB_w$ are naturally 
affine varieties (they are the set of orbits of an action of a reductive group on an affine variety with all 
orbits being of the same dimension hence closed). Similarly, by (d) above, $G^F\bsl\tX_w$ and $G^F\bsl X_w$ are 
naturally affine varieties.

The affineness properties (c),(d) can be strengthened in certain cases as follows.

\proclaim{Theorem 0.4} Assume that $G$ is almost simple of type $A_n,B_n,C_n$ or $D_n$. We assume also that
$\bul=1$. Let $w\in\WW$ be a $\bul$-elliptic element of minimal length in its $\bul$-conjugacy class. 

(a)  If $q=1$, then $G\bsl\tfB_w$ is isomorphic to $\kk^{l(w)}$ and $G\bsl\fB_w$ is isomorphic to
$T^*_w\bsl\kk^{l(w)}$ for a $T^*_w$-action on $\kk^{l(w)}$.

(b) If $q>1$, then $G^F\bsl\tX_w$ is quasi-isomorphic (see 1.1) to $\kk^{l(w)}$ and $G^F\bsl X_w$ is 
quasi-isomorphic (see 1.1) to $T^*_w\bsl\kk^{l(w)}$ for a $T^*_w$-action on $\kk^{l(w)}$.
\endproclaim
This is proved in \S4. In a sequel to this paper it is shown that (a),(b) continue to hold without the assumption
that $\bul=1$. We conjecture that (a),(b) hold for $G$ of any type.

\subhead 0.5\endsubhead
Let $w\in\WW$ and let $\d$ be the smallest integer $\ge1$ such that $\bul^\d=1$. If $q>1$,  $F:X_w@>>>X_{w^\bul}$
and $F^\d:X_w@>>>X_w$ are well defined. We propose an extension of these maps to the case of $\fB_w$ namely
$\Ps:\fB_w@>>>\fB_{w^\bul}$, $(g,B)\m(g,gBg\i)$, see 1.2; we then have $\Ps^\d:\fB_w@>>>\fB_w$. In some respects 
$\Ps,\Ps^\d$ can be viewed as analogues for $q=1$ of the Frobenius maps $F,F^\d$. Assume for example that $w$ is 
$\bul$-elliptic of minimal possible length in $\WW$. There is some evidence that, for any $i$, the 
$(i+\dim G)$-th perverse cohomology sheaf of $\r_!\bbq$ is direct sum of mutually nonisomorphic simple character 
sheaves stable under the map induced by $\Ps^\d$ and $\Ps^\d$ acts on each of these summands as multiplication by
a root of $1$ which is obtained from an eigenvalue of $F^\d$ on $H^i_c(X_w,\bbq)$ (described in \cite{\COX}) by 
$q\to1$.

\subhead 0.6\endsubhead
Let $w\in\WW$ and let $i\in I$ be such that $l(w)=l(s_iw)+1=l(s_iws_{i^\bul})$.
When $q>1$ a quasi-isomorphism (see 1.1) $\s_i:X_w@>>>X_{s_iws_{i^\bul}}$ was defined in \cite{\DL}.
In the late 1970's and early 1980's I observed (unpublished but mentioned in \cite{\MB, 5A} and \cite{\BM}) that 
by taking compositions of various $\s_i$ one can obtain nontrivial quasi-automorphisms of $X_w$ corresponding to 
elements in the stabilizer of $w$ for the $\bul$-conjugacy action (see 1.3, 1.4). Further examples of this 
phenomenon were later found by Digne and Michel \cite{\DM}. Additional examples are given in \S1,\S2. These 
examples are valid not only for $X_w$ but also for $\tX_w$, $\fB_w$ or $\tfB_w$ since in 2.3 and 2.6 we define 
quasi-isomorphisms analogous to $\s_i$ in the case when $X_w$ is replaced by $\tX_w$, $\fB_w$ or $\tfB_w$.

\subhead 0.7\endsubhead
In \S5 we give another example of the close relation between the varieties $\fB_w,X_w$ by proving (under the 
assumption that $\kk$ is as in case 2) a formula relating the number of rational points over a finite field of 
$\fB_w\T_D\fB_{w'}$ and of $G^F\bsl(X_w\T X_{w'})$. 

\subhead 0.8\endsubhead
{\it Notation.} For any $w\in\WW$ we set $\cl(w)=\{i\in I;l(s_iw)<l(w)\}$, $\car(w)=\{i\in I; l(ws_i)<l(w)\}$. 
For $k\in\ZZ$ let $w\m w^{\bul^k}$ be the $k$-th power of $\bul$. Let $w_0$ be the longest element of $\WW$. Let 
$\hat\WW$ be the braid group of $\WW$ with generators $\hs_i$ corresponding to $s_i$. If $X$ is a set and 
$f:X@>>>X$ is a map we write $X^f$ instead of $\{x\in X;f(x)=x\}$. If $X$ is finite we write $|X|$ for the 
cardinal of $X$. 

\head 1. Paths\endhead
\subhead 1.1\endsubhead
Let $C$ be a $\bul$-elliptic $\bul$-conjugacy class in $\WW$. Let $C_{min}$ be the set of elements of minimal
length of $C$. If $w\in C_{min}$ and $i\in\cl(w)$ then $w':=s_iws_i^\bul\in C_{min}$ and $i^*\in\car(w')$; we 
then write $w@>i^+>>w'$ . Conversely if $v\in C_{min}$ and $j^\bul\in\car(v)$ then $v':=s_jvs_j^\bul\in C_{min}$ 
and $j\in\cl(v')$; we then write $v@>j^->>v'$. Note that if $w,w'\in\WW$ then the conditions $w@>i^+>>w'$ and 
$w'@>i^->>w$ are equivalent. Let $\G_C$ be the graph whose vertices are the elements of $C_{min}$ and whose edges
are the triples $w\ovs i\to\smile w'$ with $w,w'$ in $C_{min}$ unordered and $i\in I$ such that either 
$w@>i^+>>w'$ or $w'@>i^+>>w$. The graph $\G_C$ has a canonical orientation in which an edge 
$w\ovs i\to\smile w'$ is oriented from $w$ to $w'$ if $w@>i^+>>w'$ and is oriented from $w'$ to $w$
if $w'@>i^+>>w$. A path in $\G_C$ is by definition a sequence $\ii$ of edges of $\G_C$ of the form
$w_1\ovs i_1\to\smile w_2\ovs i_2\to\smile,\do,\ovs i_{t-1}\to\smile w_t$. For such $\ii$ we must have 
$w_t=z_\ii\i w_1z_\ii^\bul$ where $z_\ii=s_{i_1}s_{i_2}\do s_{i_{t-1}}\in\WW$; we shall also set 
$$\ti z_\ii=\hs_{i_1}^{\e_1}\hs_{i_2}^{\e_2}\do\hs_{i_{t-1}}^{\e_{t-1}}\in\hat\WW\tag a$$
where $\e_r=1$ if $w_r@>i_r^+>>w_{r+1}$, $\e_r=-1$ if $w_r@>i_r^->>w_{r+1}$. We shall sometime specify $\ii$ by 
the symbol $[w_1;*_1,*_2,\do,*_{t-1}]$ where $*_k=i_k$ if $\e_k=1$ and $*_k=\ov i_k$ if $\e_k=-1$ ($\e_k$ as in
(a).) Note that $w_2,\do,w_t$ are uniquely determined by $w_1,i_1,i_2,\do,i_{t-1}$).

For $w,w'\in C_{min}$ let $\cp_{w,w'}$ be the set of paths in $\G_C$ such that the corresponding sequence
$w_1,w_2,\do,w_t$ satisfies $w_1=w,w_t=w'$. For example if $w=s_{i_1}s_{i_2}\do s_{i_r}$ is a reduced expression 
in $\WW$ then $[w;i_1,i_2,\do,i_r]\in\cp_{w,w^\bul}$.

The following result is due to Geck-Pfeiffer \cite{\GP, 3.2.7} (in the case where $\bul=1$) and to
Geck-Kim-Pfeiffer \cite{\GKP} and He \cite{\HE} in the remaining cases.

(b) {\it For any $w,w'\in C_{min}$, the set $\cp_{w,w'}$ is nonempty.}
\nl
For $w,w'\in C_{min}$ we identify a path $[w;*_1,*_2,\do,*_{t-1}]\in\cp_{w,w'}$ with the path
$[w;*'_1,*'_2,\do,*'_{t'-1}]\in\cp_{w,w'}$ in the following cases:

(i) $t'=t-2$, $*_k=i,*_{k+1}=\ov i$ (for some $i\in I$ and some $k$), and $*'_1,*'_2,\do,*'_{t'-1}$ is obtained 
from $*_1,*_2,\do,*_{t-1}$ by dropping $*_k,*_{k+1}$;

(ii) $t'=t-2$, $*_k=\ov i,*_{k+1}=i$ (for some $i\in I$ and some $k$), and $*'_1,*'_2,\do,*'_{t'-1}$ is obtained 
from $*_1,*_2,\do,*_{t-1}$ by dropping $*_k,*_{k+1}$;

(iii) $t'=t$, $*_k=i,*_{k+1}=j,*_{k+2}=i,\do$ ($m$ terms), $*'_k=j,*'_{k+1}=i,*'_{k+2}=j,\do$ ($m$ terms), (for 
some $i\ne j$ in $I$ with $s_is_j$ of order $m$ and some $k$) and $*'_u=*_u$ for all other indices;

(iv) $t'=t$, $*_k=\ov i,*_{k+1}=\ov j,*_{k+2}=\ov i,\do$ ($m$ terms), 
$*'_k=\ov j,*'_{k+1}=\ov i,*'_{k+2}=\ov j,\do$ ($m$ terms), (for some $i\ne j$ in $I$ with $s_is_j$ of order $m$ 
and some $k$) and $*'_u=*_u$ for all other indices.
\nl
This generates an equivalence relation on $\cp_{w,w'}$; we denote by $\bcp_{w,w'}$ the set of equivalence 
classes. For $w,w',w''\in C_{min}$, concatenation $\cp_{w,w'}\T\cp_{w',w''}@>>>\cp_{w,w''}$ induces a map
$\bcp_{w,w'}\T\bcp_{w',w''}@>>>\bcp_{w,w''}$ which makes $\sqc_{w,w'\in C_{min}}\bcp_{w,w'}$ into a groupoid. In 
particular for $w\in C_{min}$, $\bcp_{w,w}$ has a natural group structure. Now $\ii\m z_\ii$ induces a group 
homomorphism 

$\t_w:\bcp_{w,w}@>>>\WW_w:=\{z\in\WW;z\i wz^\bul=w\}$
\nl
and $\ii\m \tz_\ii$ induces a group homomorphism $\ti\t_w:\bcp_{w,w}@>>>\hat\WW$.

\subhead 1.2\endsubhead
Let $C$ be a $\bul$-elliptic $\bul$-conjugacy class in $\WW$ and let $w\in C_{min}$. We state the following
conjecture.

(a) {\it The homomorphism $\t_w:\bcp_{w,w}@>>>\WW_w$ is surjective.}
\nl
In 1.5, 1.6 we sketch a proof of (a) assuming that $\WW$ is of classical type and $\bul=1$; in 1.4 we consider in
more detail a case arising from $D_4$.

In any case, if $w^\bul=w$ then $w$ is in the image of $\t_w$. In particular, if $\WW_w$ is generated by $w$ then
(a) holds for $w$. Also from 1.1(b) we see that if (a) holds for some $w\in C_{min}$ then it holds for any
$w\in C_{min}$. We say that (a) holds for $C$ if it holds for some (or equivalently any) $w\in C_{min}$.

\subhead 1.3\endsubhead 
Assume that $w=w_0$ and $y^\bul=wyw\i$ for any $y\in\WW$. Then the $\bul$-conjugacy class of $w$ is $C=\{w\}$ and
is $\bul$-elliptic. For any $y\in\WW$ and any reduced expression $y=s_{i_1}s_{i_2}\do s_{i_k}$ for $y$, we have 
$\ii:=[y;i_1,i_2,\do,i_k]\in\cp_{w,w}$, $z_\ii=y$. Thus the image of $\t_w$ is $\WW_w=\WW$ and 1.2(a) holds in 
this case.

\subhead 1.4\endsubhead 
In the remainder of this section we assume that $\bul=1$ on $\WW$.
We will often denote an element $s_{i_1}s_{i_2}s_{i_3}\do s_{i_k}$ of $\WW$ as $i_1i_2i_3\do i_k$.

The following example appeared in the author's work (1982, unpublished).
Assume that $\WW$ is of type $D_4$. Let $S=\{s_0,s_1,s_2,s_3\}$ with $s_1,s_2,s_3$ commuting. Let $C$ be the 
conjugacy class of $\WW$ consisting of the twelve elements (of length six) $0i0j0k$ and $i0j0k0$ (where $i,j,k$ 
is a permutation of $1,2,3$). Note that $C=C_{min}$ is elliptic and any $w\in C$ has order $4$. 
We have $\cl(0i0j0k)=\{0,i\}$, $\car(0i0j0k)=\{j,k\}$, $\cl(i0j0k0)=\{i,j\}$, $\car(i0j0k0)=\{0,k\}$. 
We have $0i0j0k@>0^+>>i0j0k0$, $0i0j0k@>i^+>>0j0i0k$, $i0j0k0@>i^+>>0j0k0i$, $i0j0k0@>j^+>>i0k0j0$ for any 
$i,j,k$.

Let $w=i0j0k0\in C=C_{min}$. Now $\WW_w$ is a nonabelian group of order $16$ generated by three elements
$\a=0ij0$, $\b=jk$, $\g=i0ki0i$ satisfying 
$$\g\a\b=\a\b\g=\b\g\a=w.\tag a$$
Note that $\b$ (resp. $\a$) is the unique element of length $2$ (resp. $4$) in $\WW_w$: if $n_i$ is the number of 
elements of length $i$ in $\WW_w$ and $t$ is an indeterminate, then 
$\sum_{i\ge0}n_it^i=1+t^2+t^4+10t^6+t^8+t^{10}+t^{12}$. We have
$$\ii:=[w;\ov 0,i,j,0]\in\cp_{w,w},\ii':=[w;j,k]\in\cp_{w,w},\ii'':=[w;i,0,k,i,\ov 0,\ov i]\in\cp_{w,w},$$
and $z_\ii=\a$, $z_{\ii'}=\b$, $z_{\ii''}=\g$. 
Thus the image of $\t_w$ contains the generators $\a,\b,\g$ of $\WW_w$ hence it is equal to 
$\WW_w$ and 1.2(a) holds for $C$. Note that a relation like (a) also holds in the group $\bcp_{w,w}$:
$$\ii''\ii\ii'=\ii\ii'\ii''=\ii'\ii''\ii=[w;i,0,j,0,k,0].\tag b$$
For example,
$$\align&\ii''\ii\ii'=[w;i,0,k,i,\ov 0,\ov i,\ov 0,i,j,0,j,k]=\\&[w;i,0,k,i,\ov i,\ov 0,\ov i,i,j,0,j,k]=
[w;i,0,k,j,0,k]=[w;i,0,j,0,k,0].\endalign$$
Also $\ii,\ii',\ii''$ commute with $[w;i,0,j,0,k,0]$ in $\cp_{w,w}$.
It follows that $\tz_{\ii''},\tz_\ii,\tz_{\ii'}$ satisfy a relation like (b) in $\hat\WW$.

\subhead 1.5\endsubhead
Let $\nn$ be an integer $\ge3$. Define $n\in\NN$ by $\nn=2n$ if $\nn$ is even , $\nn=2n+1$ if $\nn$ is odd. Let 
$W$ be the group of all permutations of $[1,\nn]$ which commute with the involution $i\m\nn-i+1$ of $[1,\nn]$. 
For $i\in[1,n-1]$ define $s_i\in W$ as a product of two transpositions $i\lra i+1$, $\nn+1-i\lra\nn-i$; define 
$s_n\in W$ to be the transposition $n\lra\nn-n+1$. Then $(W,\{s_i;i\in[1,n]\})$ is a Weyl group of type $B_n$. In
this subsection we assume that $G$ is almost simple of type $C_n$ (or $B_n$) and we identify $\WW$ with $W$ with 
$\nn=2n$ (or $\nn=2n+1$) as Coxeter groups in the standard way.

Let $p_*=(p_1\ge p_2\ge\do\ge p_\s)$ be a sequence in $\ZZ_{>0}$ such that $p_1+\do+p_\s=n$. Define a partition 
$m_1+m_2+\do+m_e=\s$ by
$$p_1=p_2=\do=p_{m_1}>p_{m_1+1}=p_{m_1+2}=\do=p_{m_1+m_2}>\do.$$
For any $r\in[1,\s]$ we define a permutation $w_r$ in $W$ by
$$\align&p_{<r}+1\m p_{<r}+2\m\do\m p_{<r}+p_r\m\nn-p_{<r}-1\m\\&\nn-p_{<r}-2\m\do\m\nn-p_{<r}-p_r\m p_{<r}+1,
\endalign$$
where $p_{<r}=\sum_{r'\in[1,r-1]}p_{r'}$ and all unspecified elements are mapped to themselves. Note that 
$w_r$ is a $2p_r$-cycle and that $w_1,w_2,\do,w_\s$ are commuting with each other. Let $w=w_1w_2\do w_\s$ and let
$C$ be the conjugacy class of $w$. Note that $C$ is elliptic and $w\in C_{min}$. For every $r\in[1,\s-1]$ such 
that $p_r=p_{r+1}$ we define an involutive permutation $h_r\in W$ by
$$p_{<r}+j\m p_{<r+1}+j\m p_{<r}+j,\nn-p_{<r}-j\m\nn-p_{<r+1}-j\m\nn-p_{<r}-j\text{ for }j\in[1,p_r]$$
(all unspecified elements are mapped to themselves). Note that $h_rw_{r+1}h_r=w_r$ and $h_rw_t=w_th_r$ for all 
$t\n\{r,r+1\}$. Hence $h_rw=wh_r$. The following result is easily verified:

(a) {\it The group $\WW_w$ is generated by the elements $w_\s$, $w_r$ ($r\in[1,\s-1],p_{r+1}>p_r$) and $h_r$ 
($r\in[1,\s-1],p_r=p_{r+1}$).}
\nl
These generators satisfy the "braid group relations" of a complex reflection group of type
$$B_{m_1}^{(2p_{m_1})}\T B_{m_2}^{(2p_{m_1+m_2})}\T\do\T B_{m_e}^{(2p_{m_1+\do+m_e})}$$
(described in \cite{\MB, 3A}); the factor $B_{m_k}^{(2p_{m_1+\do+m_k})}$ is generated by $h_{m_1+\do+m_{k-1}+u}$ 
($u\in[1,m_k-1]$) and by $w_{m_1+\do+m_k}$.

It is immediate that for $r\in[1,\s]$ we have (setting $a=n-(p_\s+p_{\s-1}+\do+p_{r+1})$):
$$\ii_r:=[w;a,a+1,\do,n-1,n,n-1,\do,a-p_r+2,a-p_r+1]\in\cp_{w,w}.$$
Note that $z_{\ii_r}=w_r$.

One can verify that for $r\in[1,\s-1]$ such that $p_r=p_{r+1}=p$ we have (setting 
$a=n-(p_\s+p_{\s-1}+\do+p_{r+1})$):
$$\align&\ii'_r:=[w; a,a+1,\do,a+p-2,a-1,a,a+1\do,a+p-4,a-2,a-1,a,\do,\\&a+p-6,\do,a-p+2,a+p-1,a+p-3,\do,a-p+1,
\ov{a-p+2},\do,\ov{a+p-6},\\&\do,\ov a,\ov{a-1},\ov{a-2},\ov{a+p-4},\do,\ov{a+1},\ov a,\ov{a-1},\ov{a+p-2},\do,
\ov{a+1},\ov a]\in\cp_{w,w},\endalign$$
For example if $p=1$ we have $z_{\ii'_r}=[w;a]$; if $p=2$ we have $z_{\ii'_r}=[w;a,a+1,a-1,\ov a]$; if $p=3$ we 
have
$$z_{\ii'_r}=[w;a,a+1,a-1,a+2,a,a-2,\ov{a-1},\ov{a+1},\ov a].$$
Note that $z_{\ii'_r}=h_r$. Using (a), we see that the image of $\t_w$ contains a set of generators of $\WW_w$ 
hence 1.2(a) holds for $C$. (In the case where $p_1=p_2=\do=p_\s$, this result is due to Digne and Michel 
\cite{\DM}.) Note that any elliptic conjugacy class in $\WW$ is of the form $C$ as above. We conjecture that

(b) {\it the braid group relations satisfied by the generators in (a) remain valid as equations in
$\bcp_{w,w}$ if $w_r$ is replaced by $\ii_r$ and $h_r$ is replaced by $\ii'_r$.}
\nl
Appplying $\ti\t_w$ we would get corresponding braid group relations in $\hat\WW$ which actually can be verified.

\subhead 1.6\endsubhead
In this subsection we assume that $G$ is almost simple of type $D_n$. Let $W,s_i$ be as in 1.5 (with 
$\nn=2n\ge8$). Le $W'$ be the group of even permutations in $W$ (a subgroup of index $2$ of $W$). If $i\in[1,n-1]$
we have $s_i\in W'$ and we set $s_{(n-1)'}=s_ns_{n-1}s_n\in W'$. Then $(W',\{s_1,s_2,\do,s_{n-1},s_{(n-1)'}\})$ 
is a Weyl group of type $D_n$. We identify $\WW$ with $W'$ as Coxeter groups as in \cite{\WEU, 1.5}. Let 
$p_*=(p_1\ge p_2\ge\do\ge p_\s)$, $w_r,w,h_r$ be as in 1.5; we assume that $\s$ is even. Then $w\in W'$. Let $C'$ 
be the conjugacy class of $w$ in $W'$. Then $C'$ is elliptic and $w\in C'_{min}$. For any $r\in[1,\s]$ we have 
$w'_r:=w_rw_\s\in W'$. For any $r\in[1,\s-1]$ such that $r\in[1,\s-1],p_r=p_{r+1}$ we have $h_r\in W'$. 
If $p_{\s-1}=p_\s$ we set $h'_{\s-1}=w_\s\i h_{\s-1}w_\s$. The following result is easily verified:

(a) {\it If $p_{\s-1}>p_\s$ then $\WW_w$ is generated by the elements $w'_\s$, $w'_r$ 
($r\in[1,\s-1],p_{r+1}>p_r$) and $h_r$ ($r\in[1,\s-2],p_r=p_{r+1}$). If $p_{\s-1}=p_\s$ then $\WW_w$ is generated
by the elements $w'_\s$, $w'_r$ ($r\in[1,\s-2],p_{r+1}>p_r$), $h'_{\s-1}$ and $h_r$ ($r\in[1,\s-1],p_r=p_{r+1}$).}
\nl
These generators satisfy the "braid group relations" of a complex reflection group of type
$$B_{m_1}^{(2p_{m_1})}\T B_{m_2}^{(2p_{m_1+m_2})}\T\do\T B_{m_{e-1}}^{(2p_{m_1+\do+m_{e-1}})}\T
D_{m_e}^{(2p_{m_1+\do+m_e})}$$
(described in \cite{\MB, 3A}); the factor $B_{m_k}^{(2p_{m_1+\do+m_k})}$ (with $k<m$) is generated by
$h_{m_1+\do+m_{k-1}+u}$ ($u\in[1,m_k-1]$) and by $w'_{m_1+\do+m_k}$; if $m_e>1$ then the factor 
$D_{m_e}^{(2p_{m_1+\do+m_e})}$ is generated by $h_{m_1+\do+m_{e-1}+u}$ ($u\in[1,m_e-1]$), by $h'_{m_1+\do+m_e-1}$
and by by $w'_{m_1+\do+m_k}$; if $m_e=1$ the factor $D_{m_e}^{(2p_{m_1+\do+m_e})}$ is taken to be a cyclic group 
of order $p_\s$.For example the "braid group relation"
$$h_{\s-1}w'_\s h'_{\s-1}=w'_\s h'_{\s-1} h_{\s-1}=h'_{\s-1}h_{\s-1}w'_\s$$
holds if $m_e>1$. (Compare with 1.4(a).)

One can verify that for $r\in[1,\s]$ we have (setting $a=n-(p_\s+p_{\s-1}+\do+p_{r+1})$):
$$\align&\ii''_r=[w;a,a+1,\do,n-1,(n-1)',n-2,\do,a-p_r+2,a-p_r+1,n-1,n-2,\\&\do,n-p_\s+1]\in\cp_{w,w}.\endalign$$
Note that $z_{\ii''_r}=w'_r$. 
On the other hand for $r\in[1,\s-1],p_r=p_{r+1}$ we have $h_r=z_{\ii'_r}$ where
$\ii'_r$ is given by the same formula as in 1.5 (but viewed in $W'$); we have $\ii'_r\in\cp_{w,w}$. 
If $p_{\s-1}=p_\s=p$ then 
$$\align&\ti\ii=[w;(n-1)',n-2,n-3,\do,p+1,p,p+1,\do,n-2,p-1,p,\do,n-4,\do,\\&3,4,2,(n-1)',n-3,\do,5,3,1,
\ov 2,\ov 4,\ov 3,\do,\ov{n-4},\do,\ov{p},\ov{p-1},\ov{n-2},\do,\\&\ov{p+1},\ov{p},\ov{p+1},\do,\ov{n-3},
\ov{n-2},\ov{(n-1)'}]\in\cp_{w,w}.\endalign$$
For example if $n=10$, $p=5$ then
$$\ti\ii:=[w;9',8,7,6,5,6,7,8,4,5,6,3,4,2,9',7,5,3,1,
\ov 2,\ov 4,\ov 3,\ov 6,\ov 5,\ov 4,\ov 8,\ov 7,\ov 6,\ov 5,\ov 6,\ov 7,\ov 8,\ov{ 9'}].$$
Note that $z_{\ti\ii}=h'_{\s-1}$. Using (a), we see that the image of $\t_w$ contains a set of generators of 
$\WW_w$ hence 1.2(a) holds for $C$. Note that any elliptic conjugacy class in $\WW$ is of the form $C'$ as above.
We conjecture that

(b) {\it the braid group relations satisfied by the generators in (a) remain valid as equations in $\bcp_{w,w}$ 
if $w'_r$ is replaced by $\ii''_r$, $h_r$ is replaced by $\ii'_r$ and $h'_{\s-1}$ is replaced by $\ti\ii$.}
\nl
Appplying $\ti\t_w$ we would get corresponding braid group relations in $\hat\WW$ which actually can be verified.

\subhead 1.7\endsubhead 
In this subsection we assume that $C$ is an elliptic conjugacy class in $\WW$ such that for some $w\in C_{min}$ 
we have $w=w_1w_2\do w_r$ where $w_1,\do,w_r$ commute with each other, $l(w)=l(w_1)+l(w_2)+\do+l(w_r)$ and the 
centralizer of $w$ is generated by $w_1,\do,w_r$. (An example of this situation is the case of $w$ in 1.5 with 
$p_1>p_2>\do>p_\s$.) In this case it is immediate that $w_i$ is in the image of $\t_w$ hence 1.2(a) holds for $C$.

Another example arises for $\WW$ of type $E_8$ (with the elements of $I$ labelled as in \cite{\GP}) and with $C$ 
consisting of elements whose characteristic polynomial in the reflection representation is $(X+1)(X^7+1)$. The 
element $w=213423454234565768$ belongs to $C_{min}$ and $l(w)=18$. We have $w=s_2x=xs_2$ for some $x$ such that 
$l(x)=17$ and $s_2x^7=w_0$. (This equation holds also in $\hat\WW$.) The centralizer of $w$ is a product of a 
cyclic group of order $2$ generated by $s_2$ and a cyclic group of order $14$ generated by $x$. We see that 
1.2(a) holds for $C$.

\subhead 1.8\endsubhead 
Assume that $\WW$ is of type $E_8$. Let $C$ be the elliptic conjugacy class in $\WW$ consisting of the elements
of order $15$. We can find $w\in C_{min}$ such that $w=u^2$ where $u=12345678$ so that $l(u)=8$, $l(w)=16$. Then 
the centralizer of $w$ consists of the powers of $u$. We have $\ii:=[w;1,2,3,4,5,6,7,8]\in\cp_{w,w}$, $z_\ii=u$ 
and we see that 1.2(a) holds for $C$.

\head 2. The morphisms $\s_i,\tis_i$\endhead
\subhead 2.1\endsubhead
If $V,V'$ are algebraic varieties over $\kk$, we say that a map of sets $f:V@>>>V'$ is a quasi-morphism if:

(for $q=1$) $f$ is a morphism, or

(for $q>1$) $f$ is composition $V=V_1@>f_1>>V_2@>f_2>>\do@>f_{t-1}>>V_t=V'$ where for each $i\in[1,t-1]$, 
$f_i:V_i@>>>V_{i+1}$ is either a morphism of algebraic varieties or $V_i=V_{i+1}$ and $f_i$ is the inverse of the
Frobenius map on $V_i$ for a rational structure over a finite subfield of $\kk$.
\nl
We say that $f$ is a quasi-isomorphism if it is a quasi-morphism and has an inverse which is a quasi-morphism. If
in addition we have $V=V'$ we say that $f$ is a quasi-automorphism.

\subhead 2.2\endsubhead
Let $w\in\WW$. We define a morphism $\Ps:\XX_w@>>>\XX_{w^\bul}$ by 

$(g,B)\m(g,gBg\i)$ if $q=1$ and $B\m F(B)$ if $q>1$.
\nl
We define a morphism $\Ps:\tXX_w@>>>\tXX_{w^\bul}$ by 

$(g,g'U^*_w)\m(g,gg'd\i U^*_{w^\bul})$ (if $q=1$) and $g'U^*_w\m F(g')U^*_{w^\bul}$ (if $q>1$).
\nl
Note that each of the morphisms $\Ps$ is a quasi-isomorphism.

\subhead 2.3\endsubhead
For any $w,w',a,b\in\WW$ such that $w=ab,w'=ba^\bul$, $l(w)=l(a)+l(b)=l(w')$ we define a morphism
$\s(a):\XX_w@>>>\XX_{w'}$ by

$(g,B)\m(g,B')$ where $B'\in\cb$ is determined by the conditions $(B,B')\in\co_a$, $(B',gBg\i)\in\co_b$ (if 
$q=1$);

$B\m B'$ where $B'\in\cb$ is determined by the conditions $(B,B')\in\co_a$, \lb $(B',F(B))\in\co_b$ (if $q>1$).
\nl
(If $q>1$, the map $\s(a)$ is defined in \cite{\DL, p.107,108}.) We have a commutative diagram
$$\CD \XX_w@>\s(a)>>\XX_{w'}\\
      @V\Ps VV      @V\Ps VV \\
      \XX_{w^\bul}@>\s(a^\bul)>>\XX_{w'{}^\bul}\endCD$$
Note that for any $w\in\WW$ we have $\s(w)=\Ps:\XX_w@>>>\XX_{w^\bul}$. 

If $w,w',a,b$ are as above then $\s(b):\XX_{w'}@>>>\XX_{w^\bul}$ is defined and \lb
$\s(b)\s(a):\XX_w@>>>\XX_{w^\bul}$ 
is equal to $\Ps$. Interchanging $(a,b)$ with $(b^{\bul\i},a)$ we see that \lb
$\s(a)\s(b^{\bul\i}):X_{w'{}^{\bul\i}}@>>>X_{w'}$ is equal to $\Ps$. Thus $\s(a):\XX_w@>>>\XX_{w'}$ is a 
quasi-isomorphism.

Let $w\in\WW$ and let $i\in\cl(w)$ be such that, setting $w'=s_iws_i^\bul$, we have $l(w)=l(w')$. Then
$\s(s_i):\XX_w@>>>\XX_{w'}$ is a well defined quasi-isomorphism; we shall often write $\s_i$ instead of $\s(s_i)$.

\subhead 2.4\endsubhead
Assume that $w\in\WW$ and $i,j$ are distinct elements of $\cl(w)$. Let $m$ be the order of $s_is_j$ and 
let $v=s_is_js_i\do=s_js_is_j\do$ (both products have $m$ factors). Let $w'=vwv^\bul$ and assume that
$$l(w)=l(s_iws_i^\bul)=l(s_js_iws_i^\bul s_j^\bul)=\do=l(vwv^\bul),$$
$$l(w)=l(s_jws_j^\bul)=l(s_is_jws_j^\bul s_i^\bul)=\do=l(vwv^\bul)$$
so that the sequences of $m$ maps
$$\XX_w@>\s_i>>\XX_{s_iws_i^\bul}@>\s_j>>\XX_{s_js_iws_i^\bul s_j^\bul}@>\s_i>>\do@>>>\XX_{vwv^\bul}$$
$$\XX_w@>\s_j>>\XX_{s_jws_j^\bul}@>\s_i>>\XX_{s_is_jws_j^\bul s_i^\bul}@>\s_j>>\do@>>>\XX_{vwv^\bul}$$
are defined. We show that both compositions are equal to 
$\s_v:\XX_w@>>>\XX_{vwv^\bul}$.

Let $(g,B)\in\fB_w$ (resp. $B\in X_w$). We can find a unique sequence $B_0,B_1,\do,B_m$ in $\cb$ such that 
$B_0=B$, $(B_0,B_1)\in\co_{s_i}$, $(B_1,B_2)\in\co_{s_j}$, $(B_2,B_3)\in\co_{s_i}$, $\do$ and 
$(B_m,gBg\i)\in\co_{vw}$ (if $q=1$), $(B_m,F(B))\in\co_{vw}$ (if $q>1$). If $q=1$ we have $\s_i(g,B)=(g,B_1)$,
$\s_j(g,B_1)=(g,B_2)$, $\do$ and $\s_v(g,B)=(g,B_m)$; thus $\s_v(g,B)=\do\s_i\s_j\s_i(g,B)$ (the product has $m$ 
factors); similarly we have $\s_v(g,B)=\do\s_j\s_i\s_j(g,B)$ (the product has $m$ factors). If $q>1$ we have 
$\s_i(B)=B_1$, $\s_j(B_1)=B_2$, $\do$ and $\s_v(B)=B_m$; thus $\s_v(B)=\do\s_i\s_j\s_i(B)$ (the product has $m$ 
factors); similarly we have $\s_v(B)=\do\s_j\s_i\s_j(B)$ (the product has $m$ factors). If $q=1$, it follows 
that $\do\s_i\s_j\s_i(g,B)=\do\s_j\s_i\s_j(g,B)$ as required. If $q>1$, it follows that 
$\do\s_i\s_j\s_i(B)=\do\s_j\s_i\s_j(B)$ as required. 

\subhead 2.5\endsubhead
Assume that $w\in\WW$ and $w=s_{i_1}s_{i_2}\do s_{i_k}$ is a reduced decomposition of $w$. Let 
$$w_1=w,w_2=s_{i_2}\do s_{i_k}s_{i_1^\bul},\do, w_{k+1}=s_{i_1^\bul}\s_{i_2^\bul}\do s_{i_k^\bul}=w^\bul.$$
Assume that  $l(w_1)=l(w_2)=\do=l(w_{k+1})$. Then the sequence of maps
$$\XX_{w_1}@>\s_{i_1}>>\XX_{w_2}@>\s_{i_2}>>\do@>\s_{i_k}>>\XX_{w_{k+1}}$$
is defined. We show that the composition is equal to $\Ps:\XX_w@>>>\XX_{w^\bul}$.

Let $(g,B)\in\fB_w$ (resp. $B\in X_w$). We can find a unique sequence $B_0,B_1,\do,B_k$ in $\cb$ such that 
$B_0=B$, $(B_0,B_1)\in\co_{s_{i_1}}$, $(B_1,B_2)\in\co_{s_{i_2}}$, $\do, (B_{k-1},B_k)\in\co_{s_{i_k}}$,  and 
$B_k=gBg\i$ (if $q=1$), $B_k=F(B)$ (if $q>1$).
From the definitions we have $\s_{i_1}(g,B)=(g,B_1)$, $\s_{i_2}(g,B_1)=(g,B_2)$, $\do$, 
$\s_{i_k}(g,B_{k-1})=(g,B_k)$
if $q=1$ and  $\s_{i_1}(B)=B_1$, $\s_{i_2}(B_1)=B_2$, $\do$, $\s_{i_k}(B_{k-1})=B_k$ if $q>1$. The desired result
follows.

\subhead 2.6\endsubhead
Let $i\in I$. Let $U^*_i$ be the unique root subgroup of $U^*$ such that $(\ds^\bul)\i U^*_i\ds^\bul\not\sub U^*$
where $s=s_i$. Let $U^{*!}=\{u\in U^*;(\ds^\bul)\i u\ds^\bul\in U^*\}$. Note that any $u\in U^*$ can be 
written uniquely in the form $u=u_!u^!$ where $u_!\in U^*_i$, $u^!\in U^{*!}$ and that $u\m u_!$, $U^*@>>>U^*_i$ 
is a homomorphism. 

Now assume that $w,w',b\in\WW$ are such that $w=sb,w'=bs^\bul$, $l(w)=l(b)+1=l(w')$. Note that

(a) $\db U^*_i\db\i\sub U^*$.  
\nl
If $q=1$ we fix $g\in D$. Let $g'\in G$ be such that $g'{}\i gg'=\dw ud,u\in U^*$ (if $q=1$) and
$g'{}\i F(g')=\dw u, u\in U^*$ (if $q>1$). We set $g'_1=g'\dw u_!\db\i$. Using (a) and the definition, we have 
$$\align&g'_1{}\i gg'_1=\db u_!\i\dw\i g'{}\i gg'\dw u_!\db\i=\db u_!\i\dw\i\dw ud\dw u_!\db\i\\&=
\db u_!\i ud\ds\db u_!\db\i=\db\ds^\bul((\ds^\bul)\i u^!\ds^\bul)d(\db u_!\db\i)\in \dw'U^*dU^*=\dw'U^*d\endalign
$$
(if $q=1$),
$$\align&g'_1{}\i F(g'_1)=\db u_!\i\dw\i g'{}\i F(g')F(\dw)F(u_!)F(\db\i)=\db u_!\i\dw\i\dw uF(\dw)F(u_!)F(\db\i)
\\&=\db u^!F(\ds)F(\db u_!\db\i)=\db F(\ds)F(\ds\i)u^!F(\ds)F(\db u_!\db\i)\in\dw'U^*\endalign$$
(if $q>1$).

Now let $v\in U^*_w$. We have $v'=\dw\i v\dw\in U^*$. Using this and $w=sb,l(w)=l(b)+1$, we deduce

(b) $\ds\i v\ds\in U^*$; 
hence $(\ds^\bul)\i dvd\i\ds^\bul\in U^*$, $(dvd\i)_!=1$ (if $q=1$) and\lb
$F(\ds\i)F(v)F(\ds)\in U^*$, $F(v)_!=1$ (if $q>1$).
\nl
We have
$$(g'v)\i gg'v=v\i g'{}\i gg'v=v\i\dw udv=\dw v'{}\i udv\in\dw U^*d$$
(if $q=1$),
$$(g'v)\i F(g'v)=v\i g'{}\i F(g')F(v)=v\i\dw uF(v)=\dw v'{}\i uF(v)\in\dw U^*$$
(if $q>1$).
We define $(g'v)_1$ in terms of $g'v$ in the same way as $g'_1$ was defined in terms of $g'$. Thus we have 
$$(g'v)_1=g'v\dw v'{}\i_!u_!(dvd\i)_!\db\i=g'v\dw v'{}\i_!u_!\db\i$$
(if $q=1$),
$$(g'v)_1=g'v\dw v'{}\i_!u_!F(v)_!\db\i=g'v\dw v'{}\i_!u_!\db\i$$
(if $q>1$); we have used that $(dvd\i)_!=1$ if $q=1$ and $F(v)_!=1$ if $q>1$, see (b). 
We have $(g'v)_1=g'_1v_1$ where
$$v_1=(g'_1)\i(g'v)_1=\db u_!\i\dw\i g'{}\i g'v\dw v'{}\i_!u_!\db\i=\db u_!\i v'v'{}\i_!u_!\db\i.$$
We show that $v_1\in U^*_{w'}$. We have 
$$v_1=(\db u_!\db\i)(\ds\i v\ds)(\db v'{}\i_!u_!\db\i)$$
and this belongs to $U^*$ by (a),(b). We have 
$$\dw'{}\i v_1\dw'=(\ds^\bul)\i z\ds^\bul$$ 
where $z=u_!\i v'v'{}\i_!u_!\in U^*$. To show that $\dw'{}\i v_1\dw'\in U^*$ it is enough to observe that 
$z_!=u_!\i v'_!v'{}\i_!u_!=1$ so that $z\in U^{*!}$.

Summarizing, we see that there is a well defined morphism $\tis_i:\tXX_w@>>>\tXX_{w'}$ such that (if $q=1$) 
$(g,g'U^*_w)\m(g,g'\dw u_!\db\i U^*_{w'})$ with $u\in U^*$ given by $g'{}\i gg'=\dw ud$ and (if $q>1$) 
$g'U^*_w\m g'\dw u_!\db\i U^*_{w'}$ with $u\in U^*$ given by $g'{}\i F(g')=\dw u$. The map $\tis_i$ commutes with
the $G^F$-actions, is compatible with the $T^*_w$ and $T^*_{w'}$ actions via the isomorphism $T^*_w@>>>T^*_{w'}$,
$t\m\ds\i t\ds$ and is compatible with the map $\s_i$ (see 2.3) via the maps $\p_w$, $\p_{w'}$. In the case where
$T^*_w$ (hence $T^*_{w'}$) is finite so that $\p_w$ (resp. $\p_{w'}$) is a principal $T^*_w$- (resp. $T^*_{w'}$-)
bundle over $\XX_w$ (resp. $\XX_{w'}$) we deduce (using the fact that $\s_i:\XX_w@>>>\XX_{w'}$ is bijective) that
$\tis_i$ is bijective; it is easy to see that in this case, $\tis_i$ is a quasi-isomorphism.

\subhead 2.7\endsubhead
Assume that $c\in\WW$ and $i_1,i_2,\do,i_k\in I$ are such that each of 
$$w_1=s_{i_1}s_{i_2}\do s_{i_k}c,w_2=s_{i_2}s_{i_3}\do s_{i_k}cs_{i_1}^\bul,\do,
w_{k+1}=cs_{i_1}^\bul s_{i_2}^\bul\do s_{i_k}^\bul$$
has length $k+l(c)$.

Let $(g,g'U^*_{w_1})\in\tXX_{w_1}$ (if $q=1$), $g'U^*_{w_1}\in\tXX_{w_1}$ (if $q>1$). Using the definitions
repeatedly we see that

(a) $\tis_{i_k}\do\tis_{i_2}\tis_{i_1}(g,g'U^*_{w_1})=(g,\hg'U^*_{w_{k+1}})$ (if $q=1$)

(b) $\tis_{i_k}\do\tis_{i_2}\tis_{i_1}(g'U^*_{w_1})=\hg'U^*_{w_{k+1}}$ (if $q>1$)
\nl
where 
$$\hg'=g'\ds_{i_1}\ds_{i_2}\do\ds_{i_k}\dc \x\dc\i,$$
$$\x=u_{i_1}(\ds_{i_1}^\bul u_{i_2}(\ds_{i_1}^\bul)\i)
\do(\ds_{i_1}^\bul\ds_{i_2}^\bul\do\ds_{i_{k-1}}^\bul u_{i_k}(\ds_{i_{k-1}}^\bul)\i\do
(\ds_{i_2}^\bul)\i(\ds_{i_1}^\bul)\i),$$
with $u_{i_s}\in U^*_{i_s}$ for $s\in[1,k]$.

\subhead 2.8\endsubhead
In the setup of 2.7 we assume that $c=1$ so that $w_1=s_{i_1}s_{i_2}\do s_{i_k}$, $w_{k+1}=w_1^\bul$, $l(w_1)=k$.
We show that

(a) $\tis_{i_k}\do\tis_{i_2}\tis_{i_1}(g,g'U^*_{w_1})=\Ps(g,g'U^*_{w_1})$ (if $q=1$)

(b) $\tis_{i_k}\do\tis_{i_2}\tis_{i_1}(g'U^*_{w_1})=\Ps(g'U^*_{w_1})$ (if $q>1$).
\nl
From 2.5 we see that 

$\s_{i_k}\do\s_{i_2}\s_{i_1}(g,g'B^*g'{}\i)=(g,gg'B^*g'{}\i g\i)$ (if $q=1$)

$\s_{i_k}\do\s_{i_2}\s_{i_1}(g'B^*g'{}\i))=F(g')B^*F(g')\i$ (if $q>1$).
\nl
hence 

$\tis_{i_k}\do\tis_{i_2}\tis_{i_1}(g,g'U^*_{w_1})=(g,gg'd\i t\i U^*_{w_1^\bul})$ (if $q=1$)

$\tis_{i_k}\do\tis_{i_2}\tis_{i_1}(g'U^*_{w_1})=(F(g')t\i U^*_{w_1^\bul})$ (if $q>1$).
\nl
for some $t\in T^*_{w_1^\bul}$. Let $\hg',\x$ be as in 2.7(a),(b). We have $\hg'=g'\dw_1\dc\x\dc\i$.
If $q=1$ we have $g'{}\i gg'=\dw_1ud$ with $u\in U^*$ hence
$$gg'd\i t\i U^*_{w_1^\bul}=\hg'U^*_{w_1^\bul}=g'\dw_1\dc\x\dc\i U^*_{w_1^\bul}=gg'd\i u\i\dc\x\dc\i 
U^*_{w_1^\bul}.$$
If $q>1$ we have $g'{}\i F(g')=\dw_1u$ with $u\in U^*$ hence
$$F(g')t\i U^*_{w_1^\bul}=\hg'U^*_{w_1^\bul}=g'\dw_1\dc\x\dc\i U^*_{w_1^\bul}=F(g')u\i\dc\x\dc\i U^*_{w_1^\bul}.$$
We see that in both cases, $t\i\in u\i\dc\x\dc\i U^*_{w_1^\bul}\sub U^*$. Since $t$ is semisimple it follows that
$t=1$. This proves (a),(b).

\subhead 2.9\endsubhead
Next we assume that $w,w',i,j,m,v$ are as in 2.4. We have $w=vc,w'=cv^\bul$ where $c\in\WW$, 
$l(c)+m=l(w)=l(w')$ and the sequences of $m$ maps
$$\tXX_w@>\tis_i>>\tXX_{s_iws_i^\bul}@>\tis_j>>\tXX_{s_js_iws_i^\bul s_j^\bul}@>\tis_i>>\do@>>>\tXX_{vwv^\bul},$$
$$\tXX_w@>\tis_j>>\tXX_{s_jws_j^\bul}@>\tis_i>>\tXX_{s_is_jws_j^\bul s_i^\bul}@>\tis_j>>\do@>>>\tXX_{vwv^\bul},$$
are defined. We show:

(a) the two compositions are equal.
\nl
We apply 2.7(a),(b) with $k=m$, $i_1,i_2,\do,i_k$ equal to $i,j,i,j,\do$ and with $w_1=w,w_{k+1}=w'$. Let 
$U^*_v$ be the subgroup of $U^*$ generated by $U^*_i$ and $U^*_j$. Le $(g,g'U^*_w)\in\tXX_w$ (if $q=1$),
$g'U^*_w\in\tXX_w$ (if $q>1$). 
We have $g'{}\i gg'=\dv\dc u'u''d$ (if $q=1$), $g'{}\i F(g')=\dv\dc u'u''$ (if $q>1$) where $u'\in U^*_v$ and
$u''\in U^*\cap(\dv^\bul U^*(\dv^\bul)\i)$ are uniquely determined.

In 2.7(a),(b) we have $\hg'=g'\dv\dc\x\dc\i$ where $\x\in U^*_v$. Since
$(g,\hg'U^*_{w'})\in\tXX_{w'}$ (if $q=1$) and $\hg'U^*_{w'}\in\tXX_{w'}$ (if $q>1$), we have
$\hg'{}\i g\hg'\in\dw'U^*d$ (if $q=1$) and $\hg'{}\i F(\hg')\in\dw'U^*$ (if $q>1$).
Thus $\dc\x\i\dc\i\dv\i g'{}\i gg'\dv\dc\x\dc\i\in\dw'U^*d$ if $q=1$ and 
$\dc\x\i\dc\i\dv\i g'{}\i F(g')\dv^\bul F(\dc\x\dc\i)\in\dw'U^*$ if $q>1$.
Hence $\dc\x\i u'u''d\dv\dc\x\dc\i\in\dw'U^*d$ if $q=1$ and 
$\dc\x\i u'u''\dv^\bul F(\dc\x\dc\i)\in\dw'U^*$ if $q>1$.
We have $\dc\x\dc\i\in U^*$ hence $\x\i u'u''\in\dv^\bul U^*(\dv^\bul)\i$ in both cases.
Since $u''\in \dv^\bul U^*(\dv^\bul)\i$ we have $\x\i u'\in\dv^\bul U^*(\dv^\bul)\i$.
But we have also $\x\i u'\in U^*_v$ and $U^*_v\cap(\dv^\bul U^*(\dv^\bul)\i)=\{1\}$ hence $\x\i u'=1$ and $\x=u'$.

If we now apply 2.7(a),(b) with $k=m$, $i_1,i_2,\do,i_k$ equal to $j,i,j,j,\do$ and with $w_1=w,w_{k=1}=w'$ then 
$\hg'$ is replaced by an element  $\hg'_1=g'\dv\dc\x_1\dc\i$ where $\x_1\in U^*_v$ and by the same argument as 
above we have $\x_1=u'$. Thus $\x=\x_1$ so that $\hg'U^*_{w'}=\hg'_1U^*_{w'}$. This proves (a).

\subhead 2.10\endsubhead
Let $C$ be a $\bul$-elliptic $\bul$-conjugacy class of $\WW$ and let $w,w'\in C_{min}$. For any 
$\ii\in\cp_{w,w'}$ given by
$$w=w_1\ovs i_1\to\smile w_2\ovs i_2\to\smile\do\ovs i_{t-1}\to\smile w_t=w'$$
we define a quasi-isomorphism $T_\ii:\XX_w@>>>\XX_{w'}$ as the composition
$$\XX_{w_1}@>\s_{i_1}^{\e_1}>>\XX_{w_2}@>\s_{i_2}^{\e_2}>>\XX_{w_3}@>>>\do @>>>\XX_{w_{t-1}}
@>\s_{i_{t-1}}^{\e_{t-1}}>>\XX_{w_t}$$
and a quasi-isomorphism $\tT_\ii:\tXX_w@>>>\tXX_{w'}$ as the composition
$$\tXX_{w_1}@>\ti\s_{i_1}^{\e_1}>>\tXX_{w_2}@>\ti\s_{i_2}^{\e_2}>>\tXX_{w_3}@>>>\do@>>>\tXX_{w_{t-1}}
@>\ti\s_{i_{t-1}}^{\e_{t-1}}>>\tXX_{w_t};$$
here $\e_1,\do,\e_{t-1}$ are as in 1.1(a). Note that $T_\ii,\tT_\ii$ commute with the $G^F$-actions.
From the definitions we see that $\Ps T_\ii=T_{\ii^\bul}\Ps$ as maps $\XX_w@>>>\XX_{w'{}^\bul}$ 
and $\Ps\tT_\ii=\tT_{\ii^\bul}\Ps$ as maps $\tXX_w@>>>\tXX_{w'{}^\bul}$ where $\ii^\bul$ is given by 
$$w^\bul=w_1^\bul\ovs i_1^\bul\to\smile w_2^\bul\ovs i_2^\bul\to\smile \do\ovs i_{t-1}^\bul\to\smile w_t^\bul
=w'{}^\bul.$$
If $w\in C_{min}$ then $\ii\m T_\ii$ (resp. $\ii\m \tT_\ii$) defines a homomorphism of the group opposed to
$\bcp_{w,w}$ into the group $\cg_w$ (resp. $\ti\cg_w$) of quasi-automorphisms of $\XX_w$ (resp. of $\tXX_w$) 
which commute with the $G^F$-action. (We use 2.4, 2.9.) Hence if $q>1$ and $i\in\ZZ$ we obtain a representation 
of $\bcp_{w,w}^{opp}$ on $H^i_c(X_w,\bbq)$ and on $H^i_c(\tX_w,\bbq)$ which commutes with the $G^F$-action;
if $q=1$ and $i\in\ZZ$ we obtain a representation of $\bcp_{w,w}^{opp}$ on the $i$-th perverse cohomology 
sheaf of $\r_!\bbq$ ($\r$ as in 0.1).

\subhead 2.11\endsubhead 
Let us return to the setup of 1.4. The following relation in the group $\cg_w$ (which I found in 1982 for $X_w$) 
follows from 1.4(b):
$$T_{\ii'}T_\ii T_{\ii''}=T_{\ii''}T_{\ii'}T_{\ii}=T_{\ii}T_{\ii''}T_{\ii'}=\Ps.\tag a$$
(An analogous relation holds for $\tT_{...}$ instead of $T_{...}$.) 
Similarly, assuming that 1.5(b), 1.6(b) hold we see that in the setup of 1.5, the quasi-automorphisms
$T_{\ii_r},T_{\ii'_r}$ corresponding to the generators $w_r,h_r$ of $\WW_w$ satisfy the braid group relations 
in 1.5 and that in the setup of 1.6, the quasi-automorphisms 
$T_{\ii''_r},T_{\ii'_r},T_{\ti\ii}$ corresponding to the generators $w'_r,h_r,h'_{\s-1}$ of $\WW_w$ satisfy 
the braid group relations in 1.6. 
The apparition of braid group relations for quasi-automorphisms of $X_w$ has been predicted (in the special case 
where $w$ is regular) by Brou\'e and Michel \cite{\BM} (based on the example in 1.4, that in \cite{\RFC, p.24} 
and that for the Coxeter element in \cite{\COX}) as a part of a stronger conjecture in which the cyclotomic Hecke
algebras \cite{\MB} enter; this stronger conjecture has been verified for $C$ as in 1.5 with $p_1=p_2=\do=p_\s$ 
in \cite{\DM}.

\head 3. Proof of Theorem 0.3\endhead
\subhead 3.1\endsubhead
We prove 0.3(a). Using 1.1(b) and the quasi-isomorphisms $\tis_i$ we see that if 0.3(a) holds for some element
$\bul$-conjugate to $w$ and of the same length as 
$w$ then it will hold also for $w$. Let $\b^+$ be the braid monoid attached to the Coxeter grop $\WW$. Let 
$w_1\m\hw_1$ be the canonical imbedding $\WW@>>>\b^+$, see \cite{\GP, 4.1.1}. From the results on "good elements"
of Geck-Michel \cite{\GM}, Geck-Kim-Pfeiffer \cite{\GKP}, He \cite{\HE}, we see that, after replacing $w$ by a 
$\bul$-conjugate element of the same length as $w$, the following holds:

($*$) we can find an integer $e\ge1$ and an element $z\in\b^+$ such that 

$ww^\bul w^{\bul^2}\do w^{\bul^{e-1}}=1$ and $\hw\hw^\bul\hw^{\bul^2}\do\hw^{\bul^{e-1}}=\hat{w_0}z$ in $\b^+$. 
\nl
Thus it is enough to prove 0.3(a) for $w$ satisfying $(*)$.
Let $s_1s_2\do s_k$ be a reduced expression of $w$. Let $s'_1s'_2\do s'_f$ be a reduced expression of $w_0$. We 
can find a sequence $s''_1,s''_2,\do,s''_h$ in $S$ such that $z=\hs''_1\hs''_2\do\hs''_h$. We have
$$(\hs_1\hs_2\do\hs_k)(\hs_1^\bul\hs_2^\bul\do\hs_k^\bul)\do(\hs_1^{\bul^{e-1}}\hs_2^{\bul^{e-1}}\do
\hs_k^{\bul^{e-1}})=\hs'_1\hs'_2\do\hs'_f\hs''_1\hs''_2\do\hs''_h.$$
(The left (resp. right) hand side contains $ke$ (resp. $f+h$) elements of $S$.) We must have $ke=f+h$. Moreover by
the definition of $\b^+$ there exist $\ss^1,\ss^2,\do,\ss^m$ ($m\ge2$) such that each $\ss^r$ is a sequence 
$\ss^r_1,\ss^r_2,\do,\ss^r_{ke}$ in $S$, $\ss^1$ is the sequence 
$$s_1,s_2,\do,s_k,s_1^\bul,s_2^\bul,\do,s_k^\bul,\do,s_1^{\bul^{e-1}},s_2^{\bul^{e-1}},\do,s_k^{\bul^{e-1}},$$
($ke$ terms), $\ss^m$ is the sequence 
$$s'_1,s'_2,\do,s'_f,s''_1,s''_2,s''_h$$
and for any $r\in[1,m-1]$ the sequence $\ss^{r+1}$ is obtained from the sequence $\ss^r$ by replacing a string 
$\ss^r_{e+1},\ss^r_{e+2},\do,\ss^r_{e+u}$ of the form $s,t,s,t,\do$ ($u$ terms, $s\ne t$ in $S$, $st$ of order 
$u$ in $\WW$) by the string $t,s,t,s,\do$ ($u$ terms). 

Now let $(g,g'U^*_w)\in\tfB_w$ (if $q=1$) and $g'U^*_w\in\tX_w$ (if $q>1$). 
Let $\fZ=\{c\in G; cgc\i=g,cg'U^*_w=g'U^*_w\}$ (if $q=1$), $\fZ=\{c\in G^F;cg'U^*_w=g'U^*_w\}$ (if $q>1$).
If $c\in\fZ$ then $g'{}\i cg'\in U^*_w$ hence $c$ is unipotent. Thus $\fZ$ is a unipotent group contained in
$B:=g'B^*g'{}\i$. We define a sequence $B_0,B_1,\do,B_{ke}$ in $\cb$ by the following requirements:

$B_{ik}=g^iBg^{-i}$ (if $q=1$) and $B_{ik}=F^i(B)$ (if $q>1$) for $i\in[0,e]$, 

$(B_{ik+j-1},B_{ik+j})\in\co_{s_j^{\bul^i}}$ for $i\in[0,e-1],j\in[1,k]$.
\nl
This sequence is uniquely determined. Now conjugation by any $c\in\fZ$ preserves each of 
$B,gBg\i,g^2Bg^{-2},\do,g^eBg^{-e}$ (if $q=1$) and each of \lb 
$B,F(B),F^2(B),\do,F^e(B)$ (if $q>1$) hence (by uniqueness) it
automatically preserves each $B_v$, $v\in[0,ke]$. Thus $\fZ\sub B_v$ for any $v\in[0,ke]$.
We define a sequence $B_*^1,B_*^2,\do,B_*^m$ such that each 
$B_*^r$ is a sequence $(B_0^r,B_1^r,\do,B_{ke}^r)$ in $\cb$ satisfying $(B^r_{j-1},B^r_j)\in\co_{\ss^r_j}$ for 
$j\in[1,ke]$, as follows: $B_*^1=(B_0,B_1,\do,B_{ke})$ and for $r\in[1,m-1]$, $B_*^{r+1}$ is obtained from 
$B_*^r$ by replacing the string $B^r_a,B^r_{a+1},\do,B^r_{a+u}$ (where 
$$(\ss^r_{a+1},\ss^r_{a+2},\do,\ss^r_{a+u})=(s,t,s,t,\do)$$
as above) by the string $B^{r+1}_a,B^{r+1}_{a+1},\do,B^{r+1}_{a+u}$ defined by 
$$\align&B^{r+1}_a=B^r_a, B^{r+1}_{a+u}=B^r_{a+u}, (B^{r+1}_a,B^{r+1}_{a+1})\in\co_t, 
(B^{r+1}_{a+1},B^{r+1}_{a+2})\in\co_s,\\& (B^{r+1}_{a+2},B^{r+1}_{a+3})\in\co_t, \do.\endalign$$
(Note that $B^{r+1}_a,B^{r+1}_{a+1},\do,B^{r+1}_{a+u}$ are uniquely determined since
$(B^r_a,B^r_{a+u})\in\co_{stst\do}=\co_{tsts\do}$ and $stst\do,tsts\do$ are reduced expressions in $\WW$.) We 
note that for any $r\in[1,m]$ any Borel subgroup in the sequence $B^r_*$ is stable under conjugation by 
any $c\in\fZ$. (For $r=1$ this has been already observed. The general case follows by induction on $r$ using the 
uniqueness in the previous sentence.) In particular any Borel subgroup in the sequence $B^m_*$ is stable under 
conjugation by any $c\in\fZ$. From the definitions we see that $(B^m_0,B^m_f)\in\co_{w_0}$ that is, $B^m_0,B^m_f$
are opposed Borel subgroups. Since both are stable under conjugation by any $c\in\fZ$ we see that 
$\fZ\sub B^m_0\cap B^m_f$, a torus. Since $\fZ$ is a unipotent group we see that $\fZ=\{1\}$. This proves 0.3(a).

We prove 0.3(b). Let $(g,B)\in\fB_w$ (if $q=1$) and $B\in X_w$ (if $q>1$). 
If $q=1$ we can find $(g,g'U^*_w)\in\tfB_w$ such that $\p_w(g,g'U^*_w)=(g,B)$.
If $q>1$ we can find $g'U^*_w\in\tX_w$ such that $\p_w(g'U^*_w)=B$.
Let $\fZ_0=\{c\in G; cgc\i=g,cBc\i=B\}$ (if $q=1$), $\fZ_0=\{c\in G^F;cBc\i=B\}$ (if $q>1$).
If $c\in\fZ_0$ then $\p_w(g,cg'U^*_w)=\p_w(g,g'U^*_w)$ (for $q=1$) and $\p_w(cg'U^*_w)=\p_w(g'U^*_w)$ (for $q>1$);
hence we have $cg'U^*_w=g't\i U^*_w$ for a unique $t\in T^*_w$. Note that $c\m t$ is a group homomorphism 
$\fZ_0@>>>T^*_w$. If $c$ is in the kernel of this homomorphism then $c$ is in the isotropy group of $(g,g'U^*_w)$
(for $q=1$) and of $g'U^*_w$ (for $q>1$); hence by 0.3(a) we have $c=1$. Thus $\fZ_0@>>>T^*_w$ is injective. This 
proves (b). More precisely, we see that $\fZ_0g'U^*_w\sub g'T^*_wU^*_w$ hence $g'{}\i\fZ_0g'\sub T^*_wU^*_w$. 
Since $g'{}\i\fZ_0g'$ is a finite diagonalizable subgroup of $T^*_wU^*_w$, it is conjugate under some element of 
$U^*_w$ to a subgroup of $T^*_w$. 

We prove 0.3(c) by a method inspired by the Bonnaf\'e-Rouquier \cite{\BR} proof of 0.3(d).
We can again assume that $w$ satisfies $(*)$. Let $Y$ be the set of all sequences $(B_0,B_1,\do,B_{e-1})\in\cb^e$ such that 
$(B_i,B_{i+1})\in\co_{w^{\bul^i}}$ for $i\in[0,e-2]$. By \cite{\BR, Proposition 3}, $Y$ is an affine variety. 
Hence $D\T Y$ is an affine subvariety of $D\T\cb^e$. Let $Y'$ be the set of all 
$(g,B_0,B_1,\do,B_{e-1})\in D\T\cb^e$ such that $B_i=g^iB_0g^{-i}$ for $i\in[1,e-1]$; this is a closed subvariety
of $D\T\cb^e$. Hence $(D\T Y)\cap Y'$ is a closed subvariety of $D\T Y$ so that it is affine. The map 
$\fB_w@>>>Y'$ given by $(g,B)\m(g,B,gBg\i,g^2Bg^{-2},\do,g^{e-1}Bg^{-e+1})$ is an isomorphism of $\fB_w$ onto 
$(G\T Y)\cap Y'$. Hence $\fB_w$ is affine. Since $\tfB_w$ is a principal bundle over $\fB_w$ with (finite) group 
$T^*_w$ and $\fB_w$ is affine, we see that $\tfB_w$ is affine. This proves 0.3(c).

\proclaim{Corollary 3.2} We preserve the setup of 0.3. 

(a) If $q=1$, any isotropy group of the $U^*_w$ action $u_1:u\m\dw\i u_1\dw udu_1\i d\i$ on $U^*$ is $\{1\}$.

(b) If $q>1$, any isotropy group of the $U^*_w$ action $u_1:u\m\dw\i u_1\dw uF(u_1)\i$ on $U^*$ is $\{1\}$.
\endproclaim
We prove (a). Let $u_1\in U^*_w,u\in U^*$ be such that $\dw\i u_1\dw udu_1\i d\i=u$. We must show that $u_1=1$. 
Note that $(\dw ud,U^*_w)\in\tfB_{\dw}$ and $(u_1\dw udu_1\i,u_1U^*_w)=(\dw ud,U^*_w)$. Thus $u_1$ is in the 
isotropy group at $(\dw u,U^*_w)$ for the $G$-action on $\tfB_w$. Using 0.3(a) we deduce that $u_1=1$, as
required.

We prove (b). Let $u_1\in U^*_w,u\in U^*$ be such that $\dw\i u_1\dw uF(u_1\i)=u$. We must show that $u_1=1$. By
Lang's theorem we can find $z\in G$ such that $z\i F(z)=\dw u$. We have $u_1z\i F(z)F(u_1\i)=z\i F(z)$ that is
$zu_1z\i=F(zu_1z\i)$. We set $u'_1=zu_1z\i$ so that $u'_1\in G^F$. In the $G^F$-action on $\tX_w$, 
$u'_1\in G^F$ sends $zU^*_w\in\tX_w$ to $u'_1zU^*_w=zu_1U^*_w=zU^*_w$. Thus $u'_1$ is in the isotropy group at
$zU^*_w$ for the $G^F$-action. Using 0.3(a) we deduce that $u'_1=1$ hence $u_1=1$, as required.

\subhead 3.3\endsubhead
We preserve the setup of 0.3. Let $U^*_w\bsl\bsl U^*$ be the set of orbits of the $U^*_w$ action on $U^*$ given 
in 3.2(a) (if $q=1$) or 3.2(b) (if $q>1$). The statements (a),(b) below are immediate.

(a) If $q>1$ we have a bijection $G^F\bsl\tX_w@>\si>>U^*_w\bsl\bsl U^*$, $g'U^*_w\m\dw\i g'{}\i F(g')$ with 
inverse induced by $u\m g'U^*_w$ where $g'\in G$, $g'{}\i F(g')=\dw u$. (See \cite{\DL, 1.12}.)

(b) If $q=1$ we have a bijection $G\bsl\tfB_w@>\si>>U^*_w\bsl\bsl U^*$, $(g,g'U^*_w)\m \dw\i g'{}\i gg'd\i$ with 
inverse induced by $u\m(\dw ud,U^*)$.

\head 4. Proof of Theorem 0.4\endhead
\subhead 4.1\endsubhead
In this section we prove the assertions about $G^F\bsl\tXX_w$ in Theorem 0.4. (The assertions about 
$G^F\bsl\XX_w$ are then an
immediate consequence.) Using 3.3 we see that it is enough to consider one group in each isogeny class.
Using 0.3(a) and 1.1(b) we see that it is enough to consider a single $w$ (of minimal length) in each elliptic 
conjugacy class of $\WW$.

Let $V$ be a $\kk$-vector space of finite dimension $n\ge2$. In this subsection we assume that (if $q=1$) we have
$\hG=G=D=SL(V)$; if $q>1$ (so that $\kk$ is an algebraic closure of $F_q$) we assume that $V$ has a fixed 
$F_q$-rational structure with Frobenius map $F:V@>>>V$ (thus $V^F$ is an $n$-dimensional $F_q$-vector space) and 
that $G=SL(V)$ with the $F_q$-rational structure and Frobenius map induced by those of $V$.

Let $\o$ be a basis element of $\L^nV$ such that $F(\o)=\o$ for the map $F:\L^nV@>>>\L^nV$ given by
$v_1\we v_2\we\do\we v_n\m F(v_1)\we F(v_2)\we\do\we F(v_n)$. (Recall that if $q=1$ we have 
$F=1$.) If $q>1$ we denote by $\ss(V)$ the set of all bijective group homomorphisms $F':V@>>>V$ such that
$F'(\l v)=\l^qF'(v)$ for all $v\in V,\l\in\kk$. 

If $q>1$ let $\ss_\o(V)$ be the set of all
$F'\in\ss(V)$ such that $F'(\o)=\o$. We have $F\in\ss_\o(V)$. Note that $G$ acts on $\ss_\o(V)$ by 
$x:F'\m x F' x\i$ and that this action is transitive; the stabilizer of $F$ is $G^F$. 

Let $\cf$ be the set of all sequences $V_*=(0=V_0\sub V_1\sub V_2\sub\do\sub V_n=V)$ of subspaces of $V$ such 
that $\dim V_i=i$ for $i\in[0,n]$. Now $G$ acts naturally (transitively) on $\cf$. For any $V_*\in\cf$ we set 
$B_{V_*}=\{g\in G;gV_*=V_*\}$, a Borel subgroup of $G$. 

If $q=1$ let $Z$ be the set of all pairs $(g,V_*)\in G\T\cf$ such that 
$V_1\ne gV_1\sub V_2$, $V_2\ne gV_2\sub V_3,$ $\do,$ $V_{n-1}\ne gV_{n-1}\sub V_n$.
If $q>1$ let $Z$ be the set of all $V_*\in\cf$ such that 
$V_1\ne F(V_1)\sub V_2$, $V_2\ne F(V_2)\sub V_3,$ $\do,$ $V_{n-1}\ne F(V_{n-1})\sub V_n$.
Now $(g,V_*)\m(g,B_{V_*})$ (if $q=1$) and $V_*\m B_{V_*}$ (if $q>1$) defines an isomorphism 
$Z@>\si>>\fB_w$ (if $q=1$) or $Z@>\si>>X_w$ (if $q>1$) for a well defined Coxeter
element $w$ of length $n-1$ in $\WW$ (an elliptic element of minimal length in its conjugacy class).

If $q=1$ let $Z'$ be the set of pairs $(g,L)$ where $g\in G$ and $L$ is a line in $V$ such that
$V=\op_{i\in[0,n-1]}g^i(L)$; if $q>1$ let $Z'$ be the set of lines $L$ in $V$ such that 
$V=\op_{i\in[0,n-1]}F^i(L)$. We have an isomorphism $Z@>\si>>Z'$ given by $(g,V_*)\m(g,V_1)$ if $q=1$ and by 
$V_*\m V_1$ if $q>1$. Combining with the earlier isomorphism we obtain an isomorphism $Z'@>\si>>\fB_w$ if $q=1$
and $Z'@>\si>>X_w$ if $q>1$. (For $q>1$ the last isomorphism appears in \cite{\DL, Sec.2}.) 

If $q=1$, let $\tZ'$ be the set of pairs 
$(g,v)\in G\T V$ such that $v\we g(v)\we\do\we g^{n-1}(v)=\o$; if $q>1$, let $\tZ'$ be the set of all
$v\in V$ such that $v\we F(v)\we\do\we F^{n-1}(v)=\o$. Note that $G^F$ acts on $\tZ'$ by 
$x:(g,v)\m(xgx\i,x(v))$ (if $q=1$) 
and by $x:v\m x(v)$ (if $q>1$). Define $\p:\tZ'@>>>Z$ by $(g,v)\m(g,L)$ (if $q=1$) and by $v\m L$ (if $q>1$) 
where $L$ is the line spanned by $v$. We can identify $\tZ'$ with $\tfB_w$ (if $q=1$) or with $\tX_w$ (if $q>1$)
in a way compatible with the $G^F$-actions and so that, if $q=1$, the diagram
$$\CD \tZ'@>\si>>\tfB_w\\
       @V\p VV  @V\p_w VV  \\
        Z'@>\si>>\fB_w\endCD$$
(and the analogous diagram with $\fB_w,\tfB_w$ replaced by $X_w,\tX_w$ if $q>1$) is commutative. (For $q>1$ see
\cite{\DL, Sec.2}.) 
If $q=1$, let $\tZ''$ be the set of all $(g,v_0,v_1,\do,v_{n-1})\in G\T V^n$ such that $v_i=g^i(v_0)$ for 
$i\in[0,n-1]$, $v_0\we v_1\we\do\we v_{n-1}=\o$. If $q>1$, let $\tZ''_0$ be the set of all 
$(v_0,v_1,\do,v_{n-1})\in V^n$ such that $v_i=F^i(v_0)$ for $i\in[0,n-1]$, 
$v_0\we v_1\we\do\we v_{n-1}=\o$; let $\tZ''$ be the set of all 
$(F',v_0,v_1,\do,v_{n-1})\in\ss_\o(V)\T V^n$ such that $v_i=F'{}^i(v_0)$ for $i\in[0,n-1]$, 
$v_0\we v_1\we\do\we v_{n-1}=\o$. 

If $q=1$ we have an isomorphism $\tZ''@>\si>>\tZ'$ given by $(g,v_0,v_1,\do,v_{n-1})\m(g,v_0)$; if $q>1$ we have
an isomorphism $\tZ''_0@>\si>>\tZ'$ given by $(v_0,v_1,\do,v_{n-1})\m v_0$. Combining with the earlier 
isomorphism we obtain an isomorphism $\tZ''@>\si>>\tfB_w$ if $q=1$ and $\tZ''_0@>\si>>\tX_w$ if $q>1$. 

If $q=1$ the $G$-action on $\tfB_w$ becomes the $G$-action on $\tZ''$ given by 
$$x:(g,v_0,v_1,\do,v_{n-1})\m(xgx\i,x(v_0),x(v_1),\do,x(v_{n-1})).$$ If $q>1$ the $G^F$-action on $\tX_w$ becomes
the $G^F$-action on $\tZ''_0$ given by $x:(v_0,v_1,\do,v_{n-1})\m(x(v_0),x(v_1),\do,x(v_{n-1}))$.
If $q>1$, $G$ acts (freely) on $\tZ''$ by $x:(F',v_0,v_1,\do,v_{n-1})\m(xF'x\i,x(v_0),x(v_1),\do,x(v_{n-1}))$.
Since $G$ acts transitively on $\ss_\o(V)$ and the stabilizer of $F$ is $G^F$ we see that the space of
$G^F$-orbits on $\tZ''_0$ may be identified with the space of $G$-orbits on $\tZ''$.
We must show that the space of $G$-orbits on $\tZ''$ is an affine space for any $q$. We define 
$\tZ''@>>>\kk^{n-1}$ by $(g,v_0,v_1,\do,v_{n-1})\m(a_1,a_2,\do a_{n-1})$ if $q=1$ and by
$(F',v_0,v_1,\do,v_{n-1})\m(a_1,a_2,\do a_{n-1})$ if $q>1$ where $a_i\in\kk$ are given by 
$g^n(v_0)=a_0v_0+a_1v_1+\do+a_{n-1}v_{n-1}$ (if $q=1$) and by 
$F'{}^n(v_0)=a_0v_0+a_1v_1+\do+a_{n-1}v_{n-1}$ (if $q>1$); 
the coefficient $a_0$ is equal to $(-1)^{n-1}$. This map is constant on the orbits of $G$ hence it 
induces a map $\mu:G\bsl\tZ''@>>>\kk^{n-1}$. Next we define a map in the opposite direction 
$\t:\kk^{n-1}@>>>G\bsl\tZ''$. 

Let $(a_1,a_2,\do a_{n-1})\in\kk^n$. Let $v_0,v_1,\do,v_{n-1}$ be any basis of 
$V$ such that $v_0\we v_1\we\do\we v_{n-1}=\o$. If $q=1$ define $g\in GL(V)$ by $g(v_0)=v_1$, 
$g(v_1)=v_2$, $\do,$ $g(v_{n-2})=v_{n-1}$, $g(v_{n-1})=(-1)^{n-1}v_0+a_1v_1+\do+a_{n-1}v_{n-1}$. We have 
$(g,v_0,v_1,\do,v_{n-1})\in\tZ''$ and the $G$-orbit of this element of $\tZ''$ is independent of the choices and 
is by definition $\t(a_1,a_2,\do,a_{n-1})$. If $q>1$ we define $F'\in\ss_\o(V)$ by the requirement that 
$F'(v_0)=v_1$, $F'(v_1)=v_2$, $\do$, $F'(v_{n-2})=v_{n-1}$, 
$F'(v_{n-1})=(-1)^{n-1}v_0+a_1v_1+\do+a_{n-1}v_{n-1}$. We have
$(F',v_0,v_1,\do,v_{n-1})\in\tZ''$ and the $G$-orbit of this element of $\tZ''$ is independent of the choices and 
is by definition $\t(a_1,a_2,\do a_{n-1})$.

It is clear that $\t$ is an inverse of $\mu$. This completes the proof of Theorem 0.4 in our case.

\subhead 4.2\endsubhead
Let $V$ be a $\kk$-vector space of finite dimension $\nn\ge3$. We set $\k=0$ if $\nn$ is even, $\k=1$ if $\nn$ 
is odd and $n=(\nn-\k)/2$. Assume that $V$ has a fixed bilinear form $(,):V\T V@>>>\kk$ and a fixed quadratic 
form $Q:V@>>>\kk$ such that either

$Q=0$, $(x,x)=0$ for all $x\in V$, $V^\pe=0$;

or

$Q\ne0$, $(x,y)=Q(x+y)-Q(x)-Q(y)$ for $x,y\in V$, $Q:V^\pe@>>>\kk$ is injective.
\nl
Here, for any subspace $V'$ of $V$ we set $V'{}^\pe=\{x\in V;(x,V')=0\}$. If $Q\ne0$ it follows that $V^\pe=0$ 
unless $\k=1$ and $p=2$ in which case $\dim V^\pe=1$. If $Q=0$ we set $\e=-1$; if
$Q\ne0$ we set $\e=1$. We have $(x,y)=\e(y,x)$ for any $x,y\in V$. A subspace $V'$ of $V$ is said to be isotropic
if $(,)$ and $Q$ are zero on $V'$. In the case where $\k=0,Q\ne0$, we fix a connected component $\ci$ of the 
space of isotropic subspaces of dimension $n$ of $V$.

Let $Is(V)$ be the group of all $g\in GL(V)$ such that $(gx,gy)=(x,y)$ for all $x,y\in V$ and $Q(gx)=Q(x)$ for all 
$x\in V$ (a closed subgroup of $GL(V)$). In this section we 
assume that (if $q=1$) $G=D$ is the identity component of $Is(V)$; if $q>1$ (so that $\kk$ is an 
algebraic closure of $F_q$) we assume that $V$ has a fixed $F_q$-rational structure with Frobenius map $F:V@>>>V$
(so that $V^F$ is an $\nn$-dimensional $F_q$-vector space), that $(F(x),F(y))=(x,y)^q$ for all $x,y\in V$, that 
$Q(F(x))=Q(x)^q$ for al $x\in V$ and that $G$ is the identity component of $Is(V)$ with the $F_q$-rational 
structure and Frobenius map induced by those of $V$; in addition we assume that $G$ is $F_q$-split.

Let $\cf'$ be the set of all sequences $V_*=(0=V_0\sub V_1\sub V_2\sub\do\sub V_\nn=V)$ of subspaces of $V$ such 
that $\dim V_i=i$ for $i\in[0,\nn]$, $Q|_{V_i}=0$, $V_i^\pe=V_{\nn-i}$ for all $i\in[0,n]$ and (in the case where
$\k=0,Q\ne0$), $V_n\in\ci$. Now $G$ acts naturally (transitively) on $\cf'$.

As in 1.5, let $W$ be the group of permutations of $[1,\nn]$ which commute with the involution $i\m\nn-i+1$ of 
$[1,\nn]$. Let $V_*,V'_*$ be two sequences in $\cf'$. Let $a_{V_*,V'_*}:i\m a_i$ be the permutation of $[1,\nn]$ 
defined in \cite{\WEU, 1.4}. 
When $\k=0,Q\ne0$ let $W'$ be the group of even permutations in $W$ (a subgroup of index $2$ of $W$), see 1.6. 
Let $s_i\in W (i\in[1,n])$ be as in 1.5. Then $(W,\{s_1,s_2,\do,s_{n-1},s_n\})$ is a Weyl group of type $B_n$. 
If $\k=0,Q\ne0$, we have $s_i\in W'$ for $i\in[1,n-1]$; as in 1.6 we set $s_{(n-1)'}=s_ns_{n-1}s_n\in W'$.
 Then $(W',\{s_1,s_2,\do,s_{n-1},s_{(n-1)'}\})$ is a Weyl group of type $D_n$. 
We identify $\WW$ with $W$ (if $(1-\k)Q=0$) and with $W'$ (if $(1-\k)Q\ne0$) as Coxeter groups 
as in \cite{\WEU, 1.5}.
For any $V_*\in\cf'$ we set $B_{V_*}=\{g\in G;gV_*=V_*\}$, a Borel subgroup of $G$. 
We identify $\cf'=\cb$ via $V_*\m B_{V_*}$.

\subhead 4.3\endsubhead
In the remainder of this paper we preserve the setup of 4.2.

Let $p_*=(p_1\ge p_2\ge\do\ge p_\s)$ be a sequence in $\ZZ_{>0}$ such that $p_1+\do+p_\s=n$. In the case where 
$\k=0,Q\ne0$ we assume in addition that $\s$ is even. For any $r\in[1,\s]$ we set 
$p_{<r}=\sum_{r'\in[1,r-1]}p_{r'}$. Let $w\in W$ be the permutation of $[1,\nn]$ defined in 1.5.
If $(1-\k)Q=0$, then $w$ is elliptic in $\WW$ and it has minimal length in its conjugacy class $C$ in $\WW$.
If $\k=0,Q\ne0$, then $w\in W'=\WW$ is elliptic and it has minimal length in its conjugacy class $C'$ in $\WW$.

If $q=1$ let $Z=\{(g,V_*,V'_*)\in G\T\cf'\T\cf';V'_*=g(V_*),a_{V_*,V'_*}=w\}$. 

If $q>1$ let $Z=\{(V_*,V'_*)\in\cf'\T\cf';V'_*=F(V_*),a_{V_*,V'_*}=w\}$. 
\nl
Note that $Z=\fB_w$ (if $q=1$) and $Z=X_w$ (if $q>1$),

If $q=1$ let $\tZ'$ be the set of all sequences $(g,v_1,v_2,\do,v_\s)\in G\T V^\s$ such that

$(g^iv_t,v_r)=0$ for any $1\le t<r\le\s$, $i\in[-p_t,p_t-1]$;

$(v_r,g^iv_r)=0$ for $i\in[-p_r+1,p_r-1]$, $Q(v_r)=0$ and $(v_r,g^{p_r}v_r)=1$, $r\in[1,\s]$;

if $\k=0,Q\ne0$, the span of $g^jv_k$ $(k\in[1,\s],j\in[0,p_k-1])$ belongs to $\ci$.
\nl
(The span in the last condition is automatically an $n$-dimensional isotropic subspace.)

If $q>1$ let $\tZ'_0$ be the set of all sequences $(v_1,v_2,\do,v_\s)\in V^\s$ such that

$(F^i(v_t),v_r)=0$ for any $1\le t<r\le\s$, $i\in[-p_t,p_t-1]$;

$(v_r,F^i(v_r))=0$ for $i\in[-p_r+1,p_r-1]$, $Q(v_r)=0$ and $(v_r,F^{p_r}(v_r))=1$, $r\in[1,\s]$;

if $\k=0,Q\ne0$, the span of $F^j(v_k)$ $(k\in[1,\s],j\in[0,p_k-1])$ belongs to $\ci$.
\nl
(The span in the last condition is automatically an $n$-dimensional isotropic subspace.) Let 
$$\ct=\{(\l_1,\l_2,\do,\l_\s)\in(\kk^*)^\s;\l_r^{q^{p_r}+1}=1 \text{ for }r\in[1,\s]\},$$
a finite group isomorphic to $T^*_w$). Then if $q=1$, $\ct$ acts (freely) on $\tZ'$ by

$(\l_1,\l_2,\do,\l_\s):(g,v_1,v_2,\do,v_\s)\m(g,\l_1v_1,\l_2v_2,\do,\l_\s v_\s)$
\nl
and if $q>1$, $\ct$ acts (freely) on $\tZ'_0$ by

$(\l_1,\l_2,\do,\l_\s):(v_1,v_2,\do,v_\s)\m(\l_1v_1,\l_2v_2,\do,\l_\s v_\s)$.
\nl
Let $Z'$ (if $q=1$) and $Z'_0$ (if $q>1$) be the space of orbits of this $\ct$-action. The following result is 
equivalent to \cite{\WEU, 3.3}.

If $q=1$ we have an isomorphism $Z'@>\si>>Z$ induced by $(g,v_1,v_2,\do,v_\s)\m(g,V_*,g(V_*))$ where for any
$r\in[1,\s]$, $i\in[0,p_r]$, $V_{p_{<r}+i}$ is the subspace of $V$ spanned by $g^jv_k$ 
$(k\in[0,r-1],j\in[0,p_k-1])$ and by $g^jv_r$ $(j\in[0,i-1])$; moreover, $V_i^\pe=V_{\nn-i}$ for all $i\in[0,n]$. 

Exactly the same proof as in \cite{\WEU, 3.3} (with the action of $g$ replaced by the action of $F$) gives the
following result.

If $q>1$ we have an isomorphism $Z'_0@>\si>>Z$ induced by $(v_1,v_2,\do,v_\s)\m(V_*,g(V_*))$ where for any 
$r\in[1,\s]$, $i\in[0,p_r]$, $V_{p_{<r}+i}$ is the subspace of $V$ spanned by $F^j(v_k)$ 
$(k\in[0,r-1],j\in[0,p_k-1])$ and by $F^j(v_r)$ $(j\in[0,i-1])$; moreover, $V_i^\pe=V_{\nn-i}$ for all 
$i\in[0,n]$. 

Combining with an earlier identification we get an isomorphism $Z'@>\si>>\fB_w$ (if $q=1$) and  $Z'_0@>\si>>X_w$
(if $q>1$). Similarly we get an isomorphism $\tZ'@>\si>>\tfB_w$ (if $q=1$) and $\tZ'_0@>\si>>\tX_w$ (if $q>1$)
compatible with the isomorphism in the previous sentence and such that the $T^*_w$-action and $\ct$-action are 
compatible.

If $q>1$ let $\ss_1(V)$ be the set of all $F'\in\ss(V)$ (see 3.1) such that $(F'(x),F'(y))=(x,y)^q$ for all 
$x,y\in V$, $Q(F'(x))=Q(x)^q$ for all $x\in V$ and such that (in the case where $\k=0,Q\ne0$) $F'$ maps $\ci$ 
onto itself and (in the case where $\k=1$) $F'$ induces the same map as $F$ on $\L^\nn(V)$. Note that $G$ acts on
$\ss_1(V)$ transitively by $x:F'\m xF'x\i$ and the stabilizer of $F\in\ss_1(V)$ is $G^F$.

If $q>1$ let $\tZ'$ be the set of all sequences $(F',v_1,v_2,\do,v_\s)$ where $F'\in\ss_1(V)$ and
$v_1,v_2,\do,v_\s$ are vectors in $V$ such that

$(F'{}^i(v_t),v_r)=0$ for any $1\le t<r\le\s$, $i\in[-p_t,p_t-1]$;

$(v_r,F'{}^i(v_r))=0$ for $i\in[-p_r+1,p_r-1]$, $Q(v_r)=0$ and $(v_r,F'{}^{p_r}(v_r))=1$, $r\in[1,\s]$;

if $\k=0,Q\ne0$, the span of $F'{}^j(v_k)$ $(k\in[1,\s],j\in[0,p_k-1])$ belongs to $\ci$.
\nl
Note that $G$ acts naturally on $\tZ'$; since $G$ acts on $\ss_1(V)$ transitively and the stabilizer of
$F\in\ss_1(V)$ is $G^F$ we see that the space of $G^F$-orbits on $\tZ'_0$ can be identified with the space of
$G$-orbits on $\tZ'$. 

Let $\tZ'_1$ be the set of all collections

$(g\in G;w^r_i\in V (r\in[1,\s],i\in[0,p_r-1]); z^r_j\in V (r\in[1,\s],j\in[1,p_r])$ 

(if $q=1$),

$(F'\in\ss_1(V);w^r_i\in V (r\in[1,\s],i\in[0,p_r-1]); z^r_j\in V (r\in[1,\s],j\in[1,p_r])$ 

(if $q>1$)
\nl
such that

(a) $(w^t_i,w^r_{i'})=0$ for all $t,r,i,i'$;  

(b) $Q(w^t_i)=0$ for all $t,i$;

(c) $(z^t_j,z^r_{j'})=0$ for all $t,r,j>0,j'>0$;  

(d) $Q(z^t_j)=0$ for all $t,j>0$; 

(e) $(w^r_i,z^r_j)=(w^r_0,z^r_{i+j})^{q^i}$ if $j>0,i+j<p_r$;

(f) $(w^r_i,z^r_j)=1$  if $j>0,i+j=p_r$; 

(g) $(w^r_i,z^r_j)=0$  if $j>0,i+j>p_r$; 

(h) $(w^t_i,z^r_j)=(w^t_0,z^r_{i+j})^{q^i}$ if $j>0,i+j<p_r, t<r$; 

(i) $(w^t_i,z^r_j)=0$ if $j>0,i+j\ge p_r,t<r$;   

(j) $(w^t_i,z^r_j)=(w^t_0,z^r_{i+j})^{q^i}$ if $j>0, i+j\le p_t, t>r$;

(k) $(w^t_i,z^r_j)=0$ if $j>0,i+j>p_t,t>r$;

(l) if $\k=0,Q\ne0$, the span of $z^r_j$ $(r\in[1,\s],j\in[1,p_r])$ belongs to $\ci$;

(m) $gw^r_i=w^r_{i+1}$ for $r\in[1,\s],i\in[0,p_r-2]$, $gw^r_{p_r-1}=z^r_{p_r}$ for $r\in[1,\s]$, 
$gz^r_j=z^r_{j-1}$ for $r\in[1,\s],j\in[2,p_r]$ (if $q=1$).

(n) $F'(w^r_i)=w^r_{i+1}$ for $r\in[1,\s],i\in[0,p_r-2]$, $F'(w^r_{p_r-1})=z^r_{p_r}$ for $r\in[1,\s]$, 
$F'(z^r_j)=z^r_{j-1}$ for $r\in[1,\s],j\in[2,p_r]$ (if $q>1$).
\nl                    
If $q=1$ we have an isomorphism $\tZ'@>\si>>\tZ'_1$ given by 
$$(g,v_1,v_2,\do,v_\s)\m(g,w^r_i,z^r_j)$$
where $w^r_i=g^{-p_r+i}v_r$ for $r\in[1,\s],i\in[0,p_r-1]$, $z^r_j=g^{p_r-j}v_r$ for $r\in[1,\s],j\in[1,p_r]$.

If $q>1$ we have an isomorphism $\tZ'@>\si>>\tZ'_1$ given by 
$$(F',v_1,v_2,\do,v_\s)\m(F',w^r_i,z^r_j)$$ 
where $w^r_i=F'{}^{-p_r+i}v_r$ for $r\in[1,\s],i\in[0,p_r-1]$, $z^r_j=F'{}^{p_r-j}v_r$ for 
$r\in[1,\s],j\in[1,p_r]$.

Let $\tZ'_2$ be the set of all collections

$(w^r_i\in V (r\in[1,\s],i\in[0,p_r-1]); z^r_j\in V (r\in[1,\s],j\in[0,p_r])$ 
\nl
such that that equations (a)-(l) hold and in addition the following equations hold:

(I) $(z^t_0,w^s_{p_s-h})=(z^t_1,w^s_{p_s-h-1})^q$ for $t,s\in[1,\s]$, $h\in[1,p_s-1]$;

(II) $(z^t_0,z^s_{p_s-h})=0$ for $t,s\in[1,\s]$, $h\in[1,p_s-1]$;
             $(z^t_0,z^s_{p_s})=(z^t_1,w^s_{p_s-1})^q$ for $t,s\in[1,\s]$;

(III) $(z^t_0,z^{t'}_0)=0$ for $t<t'$ in $[1,\s]$;

(IV) $Q(z^t_0)=0$ for $t\in[1,\s]$.
\nl
It is easy to verify that the elements $w^r_i,z^r_j$ associated with a collection in $\tZ'_2$ form a basis of $V$
except if $\k=1$ when they form a basis of a hyperplane in $V$ on which $(,)$ is nondegenerate. 

We have an isomorphism $\tZ'_1@>\si>>\tZ'_2$ given by 
$$\align&(g,(w^r_i)_{r\in[1,\s],i\in[0,p_r-1]}; (z^r_j)_{r\in[1,\s],j\in[1,p_r]})\m\\&
((w^r_i)_{r\in[1,\s],i\in[0,p_r-1]}; (z^r_j)_{r\in[1,\s],j\in[0,p_r]})\endalign$$ 
(if $q=1$),
$$\align&(F',(w^r_i)_{r\in[1,\s],i\in[0,p_r-1]}; (z^r_j)_{r\in[1,\s],j\in[1,p_r]})\m\\&
((w^r_i)_{r\in[1,\s],i\in[0,p_r-1]}; (z^r_j)_{r\in[1,\s],j\in[0,p_r]})\endalign$$ 
(if $q>1$) where 

(o) $z^r_0=gz^r_1$ for $r\in[1,\s]$ (if $q=1$) and $z^r_0=F'(z^r_1)$ for $r\in[1,\s]$ (if $q>1$).
\nl 
The inverse map is given by
$$\align&((w^r_i)_{r\in[1,\s],i\in[0,p_r-1]}; (z^r_j)_{r\in[1,\s],j\in[0,p_r]})\m\\&
(g,(w^r_i)_{r\in[1,\s],i\in[0,p_r-1]}; (z^r_j)_{r\in[1,\s],j\in[1,p_r]})\endalign$$
(if $q=1$),
$$\align&((w^r_i)_{r\in[1,\s],i\in[0,p_r-1]}; (z^r_j)_{r\in[1,\s],j\in[0,p_r]})\m\\&
(F',(w^r_i)_{r\in[1,\s],i\in[0,p_r-1]}; (z^r_j)_{r\in[1,\s],j\in[1,p_r]})\endalign$$
(if $q>1$), where $g\in G$ (if $q=1$) and $F'\in\ss_1(V)$ (if $q>1$) is defined on 
$$w^r_i (r\in[1,\s],i\in[0,p_r-1],\qua z^r_j (r\in[1,\s],j\in[1,p_r])$$
by (p),(q),(o); if $\k=1$, we denote by $\x$ the unique vector in $V$ such that 
$(w^r_i,\x)=0$, $(z^r_j,\x)=0$ for all $r,i,j>0$ and such that
$$w^1_0\we w^1_{p_1-1}\we\do\we w^\s_0\we w^\s_{p_\s-1}\we z^1_1\we z^1_{p_1}\we\do\we z^\s_1\we z^\s_{p_\s}\we\x
=\o$$
(with $\o$ being a fixed basis element of $\L^\nn(V)$) and the value $g(\x)\in V$ (resp. $F'(\x)\in V$) is 
uniquely determined by the requirement that $g\in G$ (resp. $F'\in\ss_1(V)$. For future reference we note that 
$\z:=(\x,\x)$ and $\z_0=Q(\x)$ depend only on $\o$ and not on $w^r_i,z^r_j$. 

Let $\tZ'_3$ be the set of all collections
$$\align&(c^r_h\in\kk (r\in[1,\s],h\in[1,p_r-1];d^{t,r}_h\in\kk (1\le t<r\le\s;h\in[1,p_r-1]);e^{t,r}_h\in\kk \\&
(1\le r<t\le\s;h\in[1,p_t]);x^{t,r}_i\in\kk (t,r\in[1,\s],i\in[0,p_r-1]); y^{t,r}_j\in\kk \\&
(t,r\in[1,\s],j\in[1,p_r]); u^t\in\kk (t\in[1,\s]); u^t=0 \text{ unless }\k=1)\endalign$$
such that the equations $(i),(ii),(iii),(iv)$ below are satisfied.
$$\align &y^{t,s}_h+\sum_{j\in[1,h-1]}(c^s_{p_s+j-h})^{q^{p_s-h}}y^{t,s}_j
+\sum_{r,j;r<s;j\in[1,h]}(e^{s,r}_{p_s+j-h})^{q^{p_s-h}}y^{t,r}_j\\&
+\sum_{r,j;r>s;j\in[1,p_r-p_s+h-1]}(d^{s,r}_{p_s+j-h})^{q^{p_s-h}}y^{t,r}_j=M\tag i\endalign$$
for any $t,s\in[1,\s],h\in[1,p_s-1]$,
where $M=(d^{s,t}_{p_s-h})^{q^{p_s-h}}\e$ if $s<t$, $M=(e^{s,t}_{p_s-h})^{q^{p_s-h}}\e$ if $s>t$, 
$M=(c^t_{p_s-h})^{q^{p_s-h}}\e$ if $s=t$;
$$\align&x^{t,s}_h+\sum_{i\in[0,h-1]}(c^s_{p_s+i-h})^{q^i}x^{t,s}_i
+\sum_{r,i;r<s;i\in[0,h-1]}(d^{r,s}_{p_s+i-h})^{q^i}x^{t,r}_i\\&
+\sum_{r,i;r>s;i\in[0,p_r-p_s+h]}(e^{r,s}_{p_s+i-h})^{q^i}x^{t,r}_i=M'\tag ii\endalign$$   
for any $t,s\in[1,\s],h\in[0,p_s-1]$,
where $M'=(e^{s,t}_{p_s})^{q^{p_s}}\e$ if $s>t$, $h=0$, $M'=\e$ if $s=t$, $h=0$, $M'=0$ if $s<t$ or if $h>0$; 
$$\align&\sum_{r<r';i\in[0,p_r-1],j\in[1,p_{r'}];i+j<p_{r'}}(x^{t,r}_iy^{t',r'}_j+x^{t',r}_iy^{t,r'}_j\e)
(d^{r,r'}_{i+j})^{q^i}\\&
+\sum_{r>r';i\in[0,p_r-1],j\in[1,p_{r'}];i+j\le p_r}(x^{t,r}_iy^{t',r'}_j+x^{t',r}_iy^{t,r'}_j\e)
(e^{r,r'}_{i+j})^{q^i}\\&
+\sum_{r;i\in[0,p_r-1],j\in[1,p_r];i+j<p_r}(x^{t,r}_iy^{t',r}_j+x^{t',r}_iy^{t,r}_j\e)(c^r_{i+j})^{q^i}\\&
+\sum_{r;i\in[0,p_r-1],j\in[1,p_r];i+j=p_r}(x^{t,r}_iy^{t',r}_j+x^{t',r}_iy^{t,r}_j\e)+u^tu^{t'}\z=0;
\tag iii\endalign$$
for any $t<t'$ in $[1,\s]$;
$$\align&\sum_{r<r';i\in[0,p_r-1],j\in[1,p_{r'}];i+j<p_{r'}}x^{t,r}_iy^{t,r'}_j(d^{r,r'}_{i+j})^{q^i}\\&
+\sum_{r>r';i\in[0,p_r-1],j\in[1,p_{r'}];i+j\le p_r}x^{t,r}_iy^{t,r'}_j(e^{r,r'}_{i+j})^{q^i}\\&
+\sum_{r;i\in[0,p_r-1],j\in[1,p_r];i+j<p_r}x^{t,r}_iy^{t',r}_j(c^r_{i+j})^{q^i}\\&
+\sum_{r;i\in[0,p_r-1],j\in[1,p_r];i+j=p_r}x^{t,r}_iy^{t,r}_j+(u^t)^2\z_0=0\tag iv\endalign$$
for any $t\in[1,\s]$ (if $Q\ne0$).

We define $\tZ'_2@>>>\tZ'_3$ by setting

$c^r_h=(w^r_0,z^r_h)$ ($h\in[1,p_r-1]$), $d^{t,r}_h=(w^t_0,z^r_h)$ ($t<r, h\in[1,p_r-1]$), 
$e^{t,r}_h=(w^t_0,z^r_h)$ ($t>r, h\in[1,p_t]$),

and defining $x^{t,r}_i,y^{t,r}_j$ and $u^t$ (if $\k=1$) by:

$z^t_0=\sum_{r;i\in[0,p_r-1]}x^{t,r}_iw^r_i+\sum_{r;j\in[1,p_r]}y^{t,r}_jz^r_j$ (if $\k=0$)

$z^t_0=\sum_{r;i\in[0,p_r-1]}x^{t,r}_iw^r_i+\sum_{r;j\in[1,p_r]}y^{t,r}_jz^r_j+u^t\x$ (if $\k=1$)
\nl
($\x$ as in the definition of the inverse of $\tZ'_1@>>>\tZ'_2$.)
This map is well defined (the equations $(i),(ii),(iii),(iv)$ come from $I,II,III,IV$). 
Consider the fibre $\frak F$ of this map at a point of $\tZ'_3$. Then $\fF$ consists of all bases of $V$ (if 
$\k=0$) or "bases" ($=$bases with one missing element) of $V$ spanning a nondegenerate hyperplane (if $\k=1$), 
with a fixed index set, such that the value of $(,)$ at any two basis (or "basis") elements is prescribed, the
value of $Q$ at any basis (or "basis") element is prescribed, and such 
that (in the case $\k=0,Q\ne0$) the elements of type $z$ in this basis span a subspace in $\ci$. These bases 
(or "bases") clearly form a single $G$-orbit; note that the elements in such a basis (or "basis") will 
automatically satisfy the equations $I,II,III,IV$. We see that $\tZ'_3$ may be identified with the space of 
$G$-orbits on $\tZ'_2$ for the obvious (free) $G$-action.

We shall denote by $\UU$ a universal polynomial with coefficients in $\kk$ in the quantities 
$$\align&c^r_h (r\in[1,\s],h\in[1,p_r-1];d^{t,r}_h (1\le t<r\le\s;h\in[1,p_r-1]);\\&
e^{t,r}_h(1\le r<t\le\s;h\in[1,p_t]\endalign$$
and the quantities
$$y^{r,s}_{p_s} (r\le s \text{ in }[1,\s]) \text{ if }\k=0,Q=0,$$
$$y^{r,s}_{p_s} (r<s \text{ in }[1,\s]) \text{ if }\k=0,Q\ne0,$$
$$y^{r,s}_{p_s} (r<s \text{ in }[1,\s]),$$
$$u^t (t\in[1,\s]) \text{ if }\k=1.$$
We order the variables $y^{t,r}_j$ (with fixed $t$ and with $j\in[1,p_r-1]$) in the definition of $\tZ'_3$ as follows: we say that
$y^{t,r}_j<y^{t,s}_k$ if $j<k$ or $j=k,r<s$. Then in the equation $(i)$ all terms other than $y^{t,s}_h$ are 
$<y^{t,s}_h$ (for $r>s$ we have $j\le p_r-p_s+h-1\le h-1$ so that $j<h$). Therefore, from $(i)$ we see by
induction on the order above that
$$y^{t,r}_j=\UU\text{ for any }j\in[1,p_r-1].\tag p$$
We order the variables $x^{t,r}_i$ (with fixed $t$ and with $i\in[0,p_r-1]$) in the definition of $\tZ'_3$ as follows: we say that
$x^{t,r}_j<x^{t,s}_k$ if $j<k$ or $j=k,r>s$. Then in the equation $(ii)$ all terms other than $x^{t,s}_h$ are 
$<x^{t,s}_h$ (for $r>s$ we have $i\le p_r-p_s+h\le h$ so that $i\le h$). Therefore, from $(ii)$ we see by
induction on the order above that
$$x^{t,r}_i=\UU\text{ for any }i\in[0,p_r-1].\tag q$$
For $s\le t$ and $h=0$ we can write equation $(ii)$ as follows:
$$\align&x^{t,s}_0+\sum_{r;r>s;p_r=p_s}e^{r,s}_{p_s}x^{t,r}_0=\e\text{ if }t=s,\\&
x^{t,s}_0+\sum_{r;r>s;p_r=p_s}e^{r,s}_{p_s}x^{t,r}_0=0\text{ if }s<t.\tag r\endalign$$

Assuming that $Q\ne0$ we now rewrite $(iv)$  using (p),(q) (the only quantities $y^{t,s}_j$ that are not of the form 
$\UU$ are those with $j=p_s$):
$$\sum_{r>r';p_{r'}=p_r}x^{t,r}_0y^{t,r'}_{p_{r'}}e^{r,r'}_{p_{r'}}
+\sum_rx^{t,r}_0y^{t,r}_{p_r}+(u^t)^2\z_0=\UU$$
that is,
$$\sum_ry^{t,r}_{p_r}(x^{t,r}_0+\sum_{s;s>r;p_s=p_r}x^{t,s}_0e^{s,r}_{p_{r}})+(u^t)^2\z_0=\UU.$$
Using (r) this becomes
$$y^{t,t}_{p_t}+\sum_{r;r>t}y^{t,r}_{p_r}(x^{t,r}_0+\sum_{s;s>r;p_s=p_r}x^{t,s}_0e^{s,r}_{p_{r}})+(u^t)^2\z_0=\UU$$   
that is
$$y^{t,t}_{p_t}=\UU.$$
Here we have assumed that $Q\ne0$; but the same holds for $Q=0$ by the definition of $\UU$.
We now rewrite $(iii)$ for $t<t'$ using (p),(q) (again, the only quantities $y^{t,s}_j$ that are not of the form 
$\UU$ are those with $j=p_s$):
$$\sum_{r>r';p_{r'}=p_r}(x^{t,r}_0y^{t',r'}_{p_{r'}}+x^{t',r}_0y^{t,r'}_{p_{r'}}\e)e^{r,r'}_{p_{r'}}
+\sum_r(x^{t,r}_0y^{t',r}_{p_r}+x^{t',r}_0y^{t,r}_{p_r}\e)=\UU$$
that is
$$\align&\sum_ry^{t',r}_{p_r}(x^{t,r}_0+\sum_{s;s>r;p_s=p_r}x^{t,s}_0e^{s,r}_{p_{r}})\\&
+\sum_ry^{t,r}_{p_r}\e(x^{t',r}_0+\sum_{s;s>r;p_s=p_r}x^{t',s}_0e^{s,r}_{p_{r}})=\UU.\endalign$$
Using (r) this becomes
$$\align&\e y^{t',t}_{p_t}+y^{t,t'}_{p_{t'}}+
\sum_{r;r>t}y^{t',r}_{p_r}(x^{t,r}_0+\sum_{s;s>r;p_s=p_r}x^{t,s}_0e^{s,r}_{p_{r}})\\&
+\sum_{r;r>t'}y^{t,r}_{p_r}\e(x^{t',r}_0+\sum_{s;s>r;p_s=p_r}x^{t',s}_0e^{s,r}_{p_{r}})=\UU.\endalign$$
that is
$$\e y^{t',t}_{p_t}+
\sum_{r;t'>r>t}y^{t',r}_{p_r}(x^{t,r}_0+\sum_{s;s>r;p_s=p_r}x^{t,s}_0e^{s,r}_{p_{r}})=\UU.$$
This shows by induction on $t'-t$ that 
$$y^{t',t}_{p_t}\in\UU$$
for all $t<t'$. We now see that the equations defining $\tZ'_3$ are all of the form $b\in\UU$ where
$b$ is any one of the variables which do not enter in the definition of $\UU$. This shows that $\tZ'_3$
is an affine space whose dimension is equal to the number of variables which enter in the definition of 
$\UU$ that is
$$\sum_r(2r-1)p_r\text{ if }\k=0,Q=0\text{ or if }\k=1,$$
$$\sum_r(2r-1)p_r-\s\text{ if }\k=0,Q\ne0.$$
This completes the proof of Theorem 0.4.

\head 5. Counting rational points\endhead
\subhead 5.1\endsubhead
In this section we describe another example of a close relation between the varieties $\fB_w,X_w$.

Let $\ch$ be the Iwahori-Hecke algebra over $\QQ(\qq)$ ($\qq$ is an indeterminate) with basis $t_w(w\in\WW)$ and 
multiplication defined by $t_wt_{w'}=t_{ww'}$ if $w,w'\in\WW$, $l(ww')=l(w)+l(w')$ and 
$t_{s_i}^2=\qq+(\qq-1)t_{s_i}$ for $i\in I$. For any $w,w'\in\WW$ let $n_{w,w'}\in\ZZ[\qq]$ be the trace of the 
linear map $\ch\m\ch$ given by $t_y\m t_wt_{y^\bul}t_{w'{}\i}$ for all $y$.

\subhead 5.2\endsubhead
In this subsection we assume that we are in case 1 but $\kk$ is as in case 2 and we are given an 
$F_q$-rational structure on $\hat G$ with Frobenius map $\Ph:\hat G@>>>\hat G$ such that $\Ph(d)=d$ and 
$\Ph(t)=t^q$ for all $t\in T^*$. Then  $T^*,B^*,D$ are $\Ph$-stable and $\Ph$ acts trivially on $\WW$. We define 
a new $F_q$-rational structure on $G$ with Frobenius map $F:G@>>>G$ such that $F(x)=d\Ph(x)d\i$ for all $x\in G$.
Note that $G,F$ are as in case 2. Thus both $\fB_w$ and $X_w$ are well defined for $w\in\WW$. Moreover 
$w\m w^\bul$ defined in terms of $G,D$ is the same as $w\m w^\bul$ defined in terms of $G,F$. Now let $w,w'$ be 
elements of $W$. Let $\fB_w\T_D\fB_{w'}=\{((g_1,B),(g'_1,B')\in\fB_w\T\fB_{w'};g_1=g_1'\}$.

Let $G^F\bsl(X_w\T X_{w'})$ be the set of orbits of the diagonal $G^F$-action on $X_w\T X_{w'}$.
Note that for any $s\in\ZZ_{>0}$, $\Ph^s$ defines $F_{q^s}$-rational structures on $\fB_w\T_D\fB_{w'}$, 
$X_w\T X_{w'}$, $G^F\bsl(X_w\T_DX_{w'})$ with Frobenius maps denoted again by $\Ph^s$. We have the following 
result.

\proclaim{Theorem 5.3} Let $s\in\ZZ_{>0}$. Let $N_s=|(\fB_w\T_D\fB_{w'})(F_{q^s})|$ and 
$N'_s=|(G^F\bsl(X_w\T X_{w'}))(F_{q^s})|$. We have $N'_s=|G^{\Ph^s}|\i N_s=n_{w,w'}|_{\qq=q^s}$.
\endproclaim
The equality $N'_s=n_{w,w'}|_{\qq=q^s}$ is proved in \cite{\RFC, 3.8} under the additional assumption that $F^s$ 
acts trivially on $\WW$. However exactly the same proof applies without that assumption. It remains to show that
$|G^{\Ph^s}|\i N_s=n_{w,w'}|_{\qq=q^s}$. Replacing $\Ph$ by $\Ph^s$ we see that we can assume that $s=1$. Hence 
it is enough to show that $N'_1=|G^\Ph|\i N_1$. Let $G^F_\x$ be the stabilizer of $\x\in X_w\T X_{w'}$ in $G^F$. 
We have
$$\align&N'_1=|(G^F\bsl(X_w\T X_{w'}))^\Ph|
=\sum_{\x\in X_w\T X_{w'};\Ph(\x)=h\x\text{ for some }h\in G^F}|G^F_\x|/|G^F|\\&
=\sum_{\x\in X_w\T X_{w'};h\in G^F;\Ph(\x)=h\x}|G^F|\i\\&=
|G^F|\i|\{(h,B,B')\in G^F\T\cb\T\cb;(B,FB)\in\co_w,(B',FB')\in\co_{w'},\\&\Ph(B)=hBh\i,\Ph(B')=hB'h\i\}|\\&
=|G^F|\i|\{(h,B,B')\in G^F\T\cb\T\cb;(B,FB)\in\co_w,(B',FB')\in\co_{w'},\\& d\i F(B)d=hBh\i,d\i F(B')d=hB'h\i\}|.
\endalign$$
We set $h=\Ph(y)y\i$ where $y\in G$ has $|G^\Ph|$ choices. The condition $F(h)=h$ becomes 
$F(\Ph(y))F(y)\i=\Ph(y)y\i$ that is $\Ph(y\i F(y))=y\i F(y)$ (since $F\Ph=\Ph F$). We get
$$\align&N'_1=|G^F|\i|G^\Ph|\i|\{(y,B,B')\in G\T\cb\T\cb;\Ph(y\i F(y))=y\i F(y),\\&(B,F(B))\in\co_w,
(B',F(B'))\in\co_{w'},y\Ph(y\i)d\i F(B)d\Ph(y)y\i=B,\\&y\Ph(y\i)d\i F(B')d\Ph(y)y\i=B'\}|.\endalign$$
We set $B_1=y\i By$, $B'_1=y\i B'y$. We get
$$\align&N'_1=|G^F|\i|G^\Ph|\i|\{(y,B_1,B'_1)\in G\T\cb\T\cb;\Ph(y\i F(y))=y\i F(y),\\&
(yB_1y\i,F(y)F(B_1)F(y\i))\in\co_w,(yB'_1y\i,F(y)F(B'_1)F(y\i))\in\co_{w'},\\&
d\i F(B_1)d=B_1,d\i F(B'_1)d=B'_1\}|,\endalign$$
$$\align&N'_1=|G^F|\i|G^\Ph|\i|\{(y,B_1,B'_1)\in G\T\cb\T\cb;\Ph(y\i F(y))=y\i F(y),\\&
(yB_1y\i,F(y)dB_1d\i F(y\i))\in\co_w,(yB'_1y\i,F(y)dB'_1d\i F(y\i))\in\co_{w'},\\&
d\i F(B_1)d=B_1,d\i F(B'_1)d=B'_1\}|.\endalign$$
We set $z=y\i F(y)\in G^\Ph$. Note that for any $z\in G^\Ph$ there are $|G^F|$ values of $y$ satisfying 
$\Ph(y\i F(y))=y\i F(y)$. We get
$$\align&N'_1=|G^\Ph|\i|\{(z,B_1,B'_1)\in G^\Ph\T\cb\T\cb; (B_1,zdB_1d\i z\i)\in\co_w,\\&
(B'_1,zdB'_1d\i z\i)\in\co_{w'},d\i F(B_1)d=B_1,d\i F(B'_1)d=B'_1\}|.\endalign$$
We set $z'=zd\in D^\Ph$. We get
$$\align&N'_1=|G^\Ph|\i|\{(z',B_1,B'_1)\in D\T\cb\T\cb; (B_1,z'B_1z'{}\i)\in\co_w,\\&
(B'_1,z'B'_1z'{}\i)\in\co_{w'},(\Ph(z'),\Ph(B_1),\Ph(B'_1))=(z',B_1,B'_1)\}|.\endalign$$
Thus $N'_1=|G^\Ph|\i N$. The theorem is proved.

\subhead 5.4\endsubhead 
Assume in addition that $G$ is semisimple and that $w,w'$ are $\bul$-elliptic of minimal length in their 
$\bul$-conjugacy class. This guarantees that $\fB_w\T_D\fB_{w'}$ is affine and the (diagonal) $G$ action on 
$\fB_w\T_D\fB_{w'}\}$ has finite isotropy groups (see 0.3); thus all its orbits have the same dimensions so they 
are all closed and the set $G\bsl(\fB_w\T_D\fB_{w'})$ of orbits of this action is naturally an affine variety. 
Note that $\Ph$ defines an $F_q$-rational structure on $G\bsl\fB_w\T_D\fB_{w'}$. We show:

(a) {\it For any $s\in\ZZ_{>0}$, the affine varieties $G\bsl(\fB_w\T_D\fB_{w'})$, $G^F\bsl(X_w\T X_{w'})$ have 
the same number of $F_{q^s}$-rational points.}
\nl
In view of 5.3 it is enough to show that any $\Ph^s$-stable $G$-orbit on $\fB_w\T_D\fB_{w'}$ contains exactly 
$|G^{\Ph^s}|$ rational points. This follows from the fact that the isotropy group in $G$ at a point of that orbit
is finite.

\enddocument